\documentclass{article}
\pdfoutput=1
\usepackage[english]{babel}
\usepackage[percent]{overpic}
\usepackage{graphics}
\usepackage{caption}
\usepackage{adjustbox}

\usepackage{graphicx} 
\usepackage{cite}
\usepackage{amssymb}
\usepackage{amsmath}
\usepackage{amsthm}
\usepackage{graphicx}
\usepackage{appendix}
\usepackage{cases}
\usepackage{color}
\usepackage[ruled,linesnumbered]{algorithm2e}
\usepackage{algorithmic}

\usepackage[colorlinks=true, allcolors=blue]{hyperref}
\usepackage{appendix}
\usepackage{booktabs}
\usepackage{multirow}
\usepackage{cases}
\usepackage{array}

\usepackage{newtxtext,newtxmath}
\usepackage{makecell}
\usepackage{indentfirst}
\usepackage{hyperref}
\hypersetup{colorlinks = true, urlcolor = blue}
\usepackage{subfigure}

\usepackage[letterpaper,top=2cm,bottom=2cm,left=3cm,right=3cm,marginparwidth=1.75cm]{geometry}

\theoremstyle{definition}
\newtheorem{definition}{Definition}[section]

\title{Action Functional as an Early Warning Indicator in the Space of Probability Measures via Schr\"{o}dinger Bridge}

\author{
\normalsize{
Peng Zhang$^{1,2,3,}$\footnotemark[2]\ ,
Ting Gao$^{1,2,3,5,}$\footnotemark[1]\ ,
Jin Guo$^{1,2,3,}$\footnotemark[3]\ ,
Jinqiao Duan$^{4,5}$\footnotemark[4]
}\\[10pt]
\footnotesize{$^1$ School of Mathematics and Statistics, Huazhong University of Science and Technology, Wuhan, China} \\
\footnotesize{$^2$ Center for Mathematical Science, Huazhong University of Science and Technology, Wuhan, China} \\
\footnotesize{$^3$ Steklov-Wuhan Institute for Mathematical Exploration, Huazhong University of Science and Technology, Wuhan, China} \\
\footnotesize{$^4$ Department of Mathematics and Department of Physics, Great Bay University, Dongguan, China} \\
\footnotesize{$^5$ Guangdong Provincial Key Laboratory of Mathematical and Neural Dynamical Systems, Great Bay University, Dongguan, China. }}

\date{ }

\begin{document}

\maketitle

\begin{abstract}

Critical transitions and tipping phenomena between two meta-stable states in stochastic dynamical systems are a significant scientific issue. In this work, we expand the methodology of identifying the most probable transition pathway between two meta-stable states with Onsager-Machlup action functional, to investigate the evolutionary transition dynamics between two meta-stable invariant sets with Schrödinger bridge. In contrast to existing methodologies such as statistical analysis, bifurcation theory, information theory, statistical physics, topology, and graph theory for early warning indicators, we introduce a novel framework on Early Warning Signals (EWS) within the realm of probability measures that align with the entropy production rate (EPR). To validate our framework, we apply it to the Morris-Lecar model and investigate the transition dynamics between a meta-stable state and a stable invariant set (the limit cycle or homoclinic orbit) under various conditions. Additionally, we analyze real Alzheimer's data from the Alzheimer’s Disease Neuroimaging Initiative database to explore EWS indicating the transition from healthy to pre-AD states. This framework not only expands the transition pathway to encompass measures between two specified densities on invariant sets, but also demonstrates the potential of our early warning indicators for complex diseases.
\end{abstract}


\section{Introduction}

In the field of neuroscience, scientists have formulated hypotheses regarding the brain's critical dynamics, suggesting that the brain may undergo phase transitions between ordered and disordered states. Some tipping points, at which the system becomes exceptionally vulnerable, could play a crucial role in the development of neural diseases. Consequently, a pivotal inquiry emerges: how does the critical brain hypothesis intersect with pathological conditions in contrast to normal brain function?

The transition pathway stands as a pivotal instrument in elucidating critical transitions or tipping phenomena within dynamical systems subject to intricate stochastic influences \cite{duan2015introduction,gao2023stochastic}.
Strategies for the effective identification of transition pathways encompass a range of methodologies such as action functional theory \cite{wei2022optimal,chen2023detecting, guo2024deep}, path integral methods \cite{huang2021estimating} and optimal transport.
Given the inherent complexity arising from the high dimensionality of real-world data with the intricate nature of invariant manifolds in dynamical systems, the Schrödinger bridge stands as the optimal solution.

\footnotetext[2]{Email: \texttt{kazusa\_zp@hust.edu.cn}}
\footnotetext[1]{Email: \texttt{tgao0716@hust.edu.cn}}
\footnotetext[3]{Email: \texttt{jinguo0805@hust.edu.cn}}
\footnotetext[4]{Email: \texttt{duan@gbu.edu.cn}}

The Schrödinger bridge can be construed as a stochastic process aimed at identifying the optimal transition path measure between two prescribed marginal distributions. The Schrödinger bridge problem was originally introduced by Schrödinger \cite{schrodinger1932theorie} and can be construed as the optimal transport problem with entropy regularization. The study of the Schrödinger bridge problem is motivated by three main factors. Firstly, Schrödinger's original thought experiment in statistical mechanics, addressing a large deviation problem rooted in quantum mechanics, serves as a foundational impetus. Secondly, the application of Sanov's theorem \cite{sanov1961probability} and Gibbs's principle of conditioned reflection \cite{dembo2009large} facilitates the determination of the most probable Zustandverteilung (macrostate), akin to resolving a maximum entropy predicament. This scenario stands as an early and pivotal illustration of inference methodologies guiding the selection of a posterior distribution while minimizing assumptions beyond available data. The third motivation lies in perceiving these challenges as a means of regularizing the optimal mass transport (OMT) problem \cite{carlier2017convergence,marino2020optimal,leonard2012schrodinger,leonard2013survey,mikami2008optimal}, thereby addressing its computational intricacies \cite{benamou2000computational,angenent2003minimizing,gushchin2024building}. Indeed, the OMT problem represents the scenario in which the diffusion coefficient of the Schrödinger bridge tends towards zero. Subsequently, significant efforts \cite{pavon2021data,bernton2019schr} have been directed towards computing regularized OMT through algorithms of the Iterative Proportional Fitting (IPF) \cite{deming1940least} -Fortet \cite{fortet1940resolution}-Sinkhorn \cite{sinkhorn1964relationship} type.

Discovering a meaningful and computationally efficient approach to interpolate a high-dimensional probability distribution stands as a challenge of profound scientific importance, including generative modeling \cite{arjovsky2017wasserstein}, transfer learning \cite{neyshabur2020being} and cellular biology \cite{schiebinger2019optimal}. Currently, an increasing body of research leverages the Schrödinger bridge methodology to address this challenge.  Pavon et al. \cite{pavon2021data} introduce an iterative procedure that employs constrained maximum likelihood estimation instead of nonlinear boundary couplings. They utilize importance sampling to propagate the functions $\phi$ and $\psi$ in solving the Schrödinger system, drawing inspiration from IPF-Fortet-Sinkhorn type algorithms. Wang et al. \cite{wang2021deep} present a two-stage generative model employing entropy interpolation with the Schrödinger bridge. They develop their algorithm by incorporating the drift term, estimated by a deep score estimator and a deep density ratio estimator, into the Euler-Maruyama method. Bortoli et al. \cite{de2021diffusion} introduce a novel framework that provides an approximation of the IPF method for addressing the Schrödinger bridge problem. This framework presents a computational tool for optimal transport, serving as a continuous state-space counterpart to the well-known Sinkhorn algorithm. Chen et al. \cite{chen2023schrodinger} establish a computationally feasible Schrödinger bridge connecting paired data by introducing a versatile reference stochastic differential equation. This approach enables the generation of mel-spectrograms from latent representations of clean text. Chen et al. \cite{chen2021likelihood} present a computational framework for likelihood training of Schrödinger bridge models, which is the most similar method to ours. Drawing from the theory of forward-backward stochastic differential equations, they transform the optimality conditions of the Schrödinger bridge into a system of stochastic differential equations, facilitating efficient model training. By leveraging the Schrödinger bridge theory, these methods enhance the efficiency of computing transition measures between distributions to a certain degree.


Various dynamical systems, spanning ecosystems \cite{tylianakis2014tipping}, climate \cite{ashwin2012tipping, kerr2015climate}, economic structures \cite{gualdi2015tipping} and social frameworks \cite{dakos2014critical, dai2012generic} can manifest a tipping point, triggering a complete system collapse \cite{gladwell2006tipping,scheffer2009early}. Forecasting these tipping points represents a significant and exceedingly arduous challenge. Recent examples garnering significant interest involve neurodegenerative disorders, such as Alzheimer's disease (AD) \cite{simons2023tipping}. 
Numerous studies \cite{jagust2018imaging,ridha2006tracking} have been undertaken to identify the change point marking the onset of accelerated atrophy in AD.  These investigations have revealed that the rate of brain atrophy in AD does not exhibit a linear progression but rather undergoes a gradual acceleration several years prior to the manifestation of symptoms \cite{kinnunen2018presymptomatic}. Simons et al. \cite{simons2023tipping} propose that tipping points and their associated elements could offer a valuable concept for enhancing the comprehension of  AD. Tipping points, align with the concept that minor alterations at a specific juncture can yield substantial and enduring repercussions within a system, encapsulating the idea that small changes can yield significant outcomes. The detection of early warning signals for imminent critical transitions is crucial, given the substantial challenges involved in reverting a system to its previous state following such events\cite{folke2004regime,scheffer2001catastrophic}. The phenomenon of critical slowing down underpins various effective indicators in Early Warning Signals (EWS), encompassing heightened autocorrelation \cite{scheffer2009early}, variance \cite{maturana2020critical, boers2018early}, and higher-order moments such as skewness and kurtosis \cite{wright2011problematic} of the system's response amplitude distribution, which fluctuates as the system approaches a critical point.  These metrics are derived from data employing a sliding window technique. The relevant quantiles are assessed within a window comprising multiple data points in the time series, sliding forward at a suitable time increment to capture temporal variations. Certain scenarios exist where transitions occur without preceding critical slowing down \cite{boerlijst2013catastrophic,boettiger2013early}. Subsequent studies tackle this constraint by employing multivariate data analysis to assess false alarms \cite{streeter2013anticipating}, incorporating metrics such as link density \cite{yang2022critical}, clustering coefficient \cite{van2013interaction}, characteristic path length \cite{godavarthi2017recurrence}, etc. Other studies detect early warning signals through dynamic changes of the network structure\cite{liu2021predicting, zhong2022identifying} and from the inconsistency of predictions and the fluctuation of hidden state representations\cite{tong2023earthquake, chen2020autoreservoir }.In contrast to these methods, we aim to identify suitable indicators within the realm of probability measures, as certain rare events are more readily detectable within the framework of probability measure space.

In summary, the goal of this study is to model and analyze transition dynamics, including brain activity, to detect early warnings of abrupt changes within these systems. Our main contributions are
\begin{itemize}
    \item Early Warning Indicator: Define an action functional, i.e. the Wasserstein distance in Benamou-Brenier representation, as an early warning indicator in probability measure space. 
    \item Transition Pathway between Invariant Manifolds: Develop an efficient way to detect transition paths between two compact invariant manifolds. 
    \item Validation: Validate our method on real Alzheimer's Disease Neuroimaging Initiative (ADNI) data and the Morris-Lecar neuron model which describes a variety of oscillatory voltage patterns of Barnacle muscle fibers under different current.
\end{itemize}


The rest of this paper is organized as follows: Section 2 provides a review of the foundational concepts of the Schrödinger bridge and its equivalent formulations. Subsequently, we introduce a novel tipping indicator in probability space, derived from the action functional of the path measure of the Schrödinger bridge. In Section 4, we validate our framework through numerical experiments using real ADNI data. Lastly, our conclusion is summarized in Section 5.

\section{Background on Schrödinger bridge }

The most probable stochastic evolution between two prescribed marginal distributions $\rho_0$ and $\rho_1$ corresponds to the solution of the Schrödinger bridge problem:
\begin{equation}
    \mathcal{P}_{SB} := \mbox{argmin} \{D_{KL}(P\|W)|P\in \mathcal{D}(\rho_0,\rho_1)\},
\end{equation}
where $\mathcal{D}$ is the space of probability measures on $\Omega = C([0,1],\mathcal{X})$, and $W\in M_+({\Omega})$ is the reference path measure on $\Omega$. Here $\mathcal{X}$ is a complete connected Riemannian manifold without boundary, $M_+({\Omega})$ denotes positive measures on $\Omega$.

Consider the Schrödinger bridge problem with prior reference measure $\overline{W}$, which has the canonical coordinate process as the solution of the following stochastic differential equation
\begin{equation}
    \mbox{d}X_t = f(t, X_t)\mbox{d}t + \sigma_t \mbox{d}\overline{W}_t, \label{prior sde}
\end{equation}
where $f$ is the drift term and $\sigma_t$ denotes the noisy intensity. 
The optimal control formulation of this problem can be written as \cite{chen2021stochastic}
\begin{numcases}{}
    &$\min\limits_{u\in \mathcal{U}} \mathcal{J}(u)=\mathbb{E}\Big[\int_0^1 \frac{1}{2\sigma_t^2}\|u_t\|^2\mbox{d}t\Big],$\\
    &$\mbox{d}X_t = [f(t,X_t)+u_t]\mbox{d}t+\sigma_t \mbox{d}\overline{W}_t$, \\
    &$X_0\sim \rho_0(x), X_1\sim \rho_1(y),$
\end{numcases}
where the family $\mathcal{U}$ consists of adapted, finite-energy control functions.

Now, we present several equivalent versions for the Schrödinger bridge problem.

\subsection{Coupled partial differential equations with boundary conditions}
Let $(\varphi,\hat{\varphi})$ be the solution of the following coupled partial differential equations (PDE)
\begin{align}
    &\frac{\partial \varphi_t}{\partial t}+ f\cdot \nabla \varphi_t + \frac{\sigma_t^2}{2}\Delta \varphi_t=0,\\
    &\frac{\partial \hat{\varphi}_t}{\partial t}+\nabla\cdot (f\hat{\varphi}_t)-\frac{\sigma_t^2}{2}\Delta \hat{\varphi}_t=0
\end{align}
with boundary conditions
\begin{align}
    &\varphi(0, x) \cdot \hat{\varphi}(0, x)=\rho_0(x),\\
    &\varphi(1, y) \cdot \hat{\varphi}(1, y)=\rho_1(y).
\end{align}
Moreover, the solution can be expressed by the path measure of the following forward stochastic differential equations (SDE) 
\begin{equation}
    \mbox{d}X_t=\big(f(t,X_t)+\sigma_t^2\nabla_x \ln{\varphi(t,X_t} )\big)\mbox{d}t + \sigma_t \mbox{d}W_t,\,\, X_0\sim \rho_0,
\end{equation}
and the backward SDE
\begin{equation}
    \mbox{d}X_t=\big(f(t,X_t)-\sigma_t^2\nabla_x \ln{\hat{\varphi}(t,X_t} )\big)\mbox{d}t + \sigma_t \mbox{d}W_t,\,\, X_1\sim \rho_1.
\end{equation}

\subsection{Relationship between Benamou-Brenier problem and Schrödinger bridge}
The dynamic version of the OMT problem was already accomplished by Benamou and Brenier \cite{benamou2000computational} by 
\begin{equation}
    \begin{aligned}
& \mathcal{W}_2\left(\rho_0, \rho_1\right)=\inf _{(\mu, v)} \int_0^1 \int_{\mathbb{R}^n}\|v(t, x)\|^2 \mu_t(d x) \mbox{d} t, \\
& \frac{\partial \mu_t}{\partial t}+\nabla \cdot(v \mu_t)=0, \\
& \mu_0=\rho_0, \quad \mu_1=\rho_1 .\label{bb}
\end{aligned}
\end{equation}

The fluid-dynamic formulation for the Schrödinger bridge can be derived as in \cite{chen2016relation} 
\begin{equation}
    \begin{aligned}
& \inf _{(\rho, v)} \int_{\mathbb{R}^n} \int_0^1\left[\frac{1}{2}\|v(t, x)\|^2+\frac{\sigma_t^2}{8}\|\nabla \log \rho\|^2\right] \rho(t, x) \mbox{d} t \mbox{d} x \\
& \frac{\partial \rho}{\partial t}+\nabla \cdot(v \rho)=0, \\
& \rho(0, x)=\rho_0(x), \quad \rho(1, y)=\rho_1(y) .\label{fluid}
\end{aligned}
\end{equation}
The fluid-dynamic formulation (\ref{fluid}) could be viewed as a regularization of the Benamou-Brenier problem (\ref{bb}). As $\sigma_t \rightarrow 0$, the solution to the problem (\ref{fluid}) converges to the solution of the Benamou-Brenier problem (\ref{bb}) as shown in \cite{benamou2000computational}.
By variable substitution, the fluid-dynamic formulation (\ref{fluid}) is equivalent to 
\begin{equation}
    \begin{aligned}
& \inf _{(\rho, v)} \int_{\mathbb{R}^n} \int_0^1\left[\frac{1}{2}\|v(t, x)\|^2\right] \rho(t, x) \mbox{d} t \mbox{d} x \\
& \frac{\partial \rho}{\partial t}+\nabla \cdot(v \rho)-\frac{\sigma^2_t}{2}\Delta \rho=0, \\
& \rho(0, x)=\rho_0(x), \quad \rho(1, y)=\rho_1(y) .\label{fluid furmulation}
\end{aligned}
\end{equation}

\subsection{Forward and backward generators}
Denoting $ R\in M_+(\Omega)$ as a Markov measure, its forward stochastic derivative $\partial+\overrightarrow{L}^W$ is defined by
\begin{equation}
    \left[\partial_t+\overrightarrow{L}_t^W\right](u)(t, x):=\lim _{h \downarrow 0} h^{-1} \mathbb{E}_W\left(u\left(t+h, X_{t+h}\right)-u\left(t, X_t\right) \mid X_t=x\right)
\end{equation}
for every measurable function $u:[0,1] \times \mathcal{X} \rightarrow \mathbb{R}$ in the set dom $\overrightarrow{L}^W$ for which the limit
exists $W_t -a.e.$ for all $0 \leq t<1$. Since the time reversed $W^*$ of $W$ is still Markov, $W$ admits a backward stochastic derivative $-\partial +\overleftarrow{L}^W$ which is defined by
\begin{equation}
    \left[-\partial_t+\overleftarrow{L}_t^W\right] u(t, x):=\lim _{h \downarrow 0} h^{-1} \mathbb{E}_W\left(u\left(t-h, X_{t-h}\right)-u\left(t, X_t\right) \mid X_t=x\right)
\end{equation}
for every measurable function $u:[0,1] \times \mathcal{X} \leftarrow \mathbb{W}$ in the set dom $\overleftarrow{L}^W$ for which the limit
exists $W_t -a.e.$ for all $0 < t\leq 1$.

Denote $\overrightarrow{L}^W=\overleftarrow{L}^W=L$ without the superscript $W$ and without the time arrows, since $W$ is assumed to be reversible.

To give the expressions of the generators $\overrightarrow{A}$
and $\overleftarrow{A}$, we should introduce the carr\'{e} du champ of $W$ as
\begin{equation}
    \Gamma(u, v):=L(u v)-u L v-v L u.
\end{equation}
In general, the forward and backward generators $\left(\partial_t+\overrightarrow{A}_t\right)_{0 \leq t \leq 1}$ and $\left(-\partial_t+\overleftarrow{A}_t\right)_{0 \leq t \leq 1}$ of
path measure $P$ depend explicitly on $t$.

Under some hypotheses in \cite{leonard2013survey} on measure $W$ , the forward and backward generators of $P$ are given for any function $u:[0,1] \times \mathcal{X} \rightarrow \mathbb{R}$ belonging to some class $\mathcal{U}$ of regular functions, by
\begin{equation}
    \begin{cases}\overrightarrow{A}_t u(x)=L u(x)+\frac{\Gamma\left(\varphi_t, u\right)(x)}{\varphi_t(x)}, & (t, x) \in[0,1) \times \mathcal{X} \\ \overleftarrow{A}_t u(x)=L u(x)+\frac{\Gamma\left(\hat{\varphi}_t, u\right)(x)}{\hat{\varphi}_t(x)}, & (t, x) \in(0,1] \times \mathcal{X}\end{cases}
\end{equation}
where $\varphi_t, \hat{\varphi}_t$ are defined by
\begin{equation*}
    \left\{\begin{array}{l}
\varphi_t(z):=\mathbb{E}_W\left(\varphi_0\left(X_0\right) \mid X_t=z\right) \\
\hat{\varphi}_t(z):=\mathbb{E}_W\left(\hat{\varphi}_1\left(X_1\right) \mid X_t=z\right)
\end{array} \quad, \quad \text { for } P_t \text {-a.e. } z \in \mathcal{X}\right. \text {. }
\end{equation*}
Rigorous statement and proof are given in \cite{léonard2011stochastic}. 

\subsection{Forward-backward stochastic differential equations}
By nonlinear Feynman-Kac formulation \cite{exarchos2018stochastic}, the Schrödinger bridge problem can be derived as the following forward-backward SDEs representation \cite{chen2021likelihood}
\begin{numcases}{}
   & $\mbox{d}X_t = (f + \sigma_t^2 \nabla \ln{\varphi_t}) \mbox{d}t + \sigma_t \mbox{d}W_t$ \label{forw} \\
& $\mbox{d}Y_t= \frac{1}{2}\sigma_t^2 (\nabla \ln{\varphi_t})^2 \mbox{d}t + \sigma_t\nabla \ln{\varphi_t} \mbox{d}W_t$ \label{backw 1} \\
& $\mbox{d} \widehat{Y}_t =(\frac{1}{2}\sigma_t^2(\nabla \ln{\widehat{\varphi}_t})^2+\nabla_x\cdot (\sigma_t^2 \nabla\ln{\widehat{\varphi}_t}-f)+ \sigma_t^2\nabla \ln{\widehat{\varphi}_t}\nabla \ln{\varphi_t})\mbox{d}t + \sigma_t\nabla \ln{\widehat{\varphi}_t} \mbox{d}W_t$ \label{backw 2}
\end{numcases}
with boundary conditions $X_0=x_0, x_0\sim \rho_0$ and $Y_1 + \widehat{Y}_1 = \ln{\rho_1(X_1)}$. Here (\ref{forw}) and
(\ref{backw 1})-(\ref{backw 2}) 
respectively represent the forward and backward SDEs. 

In the following example, we employ Chen et al. \cite{chen2021likelihood} as our computational framework for finding the optimal transition path measure in the forward-backward SDEs. The loss function is defined by the log-likelihood of the Schrödinger bridge problem as
\begin{equation}
    \mathcal{L} = - \ln{p_0(x_0)}= \int_0^1 \mathbb{E}\big[\frac{1}{2}(\|\sigma_t \nabla\varphi_t\| +\|\sigma_t \nabla\widehat{\varphi}_t\|)^2 + \nabla_x\cdot (\sigma_t^2 \nabla\widehat{\varphi}_t-f)\big]\mbox{d}t - \mathbb{E}[\ln{p_1(X_1)}], \label{loss function}
\end{equation}
where the expectation is taken over the forward SDE (\ref{forw}) with the initial condition $X_0 = x_0$. The training details are shown in Algorithm \ref{algorithm SB}.

\begin{algorithm}[t]
\caption{Likelihood training of SB-FBSDE}\label{algorithm SB}
\KwIn{boundary distributions $p_{\mbox{data}}$ and $p_{\mbox{prior}}$, parameterized policies $Z(\cdot,\cdot;\theta)$ and $\hat{Z}(\cdot,\cdot;\phi)$}
    \For{$k = 1$ \KwTo{$K$}}
    {   
     Sample $X_{t\in [0,T]}$ from (\ref{forw}), where $x_0\sim p_{\mbox{data}}$ (computational graph retained)\;
     compute $\mathcal{L}$ with \eqref{loss function}\;
     Update($\theta, \phi$) with $\nabla_{\theta,\phi}\mathcal{L}_{SB}(x_0;\theta,\phi)$
    }
\end{algorithm}

\subsection{Equivalence between the Onsager-Machlup action functional and Schrödinger bridge}
One of the ways to characterize the most probable transition pathway of a stochastic process between two metastable states is the Onsager-Machlup action functional. This scenario can be viewed as a specific instance of identifying the most probable transition path between two specified marginal distributions, i.e. two Dirac $\delta$ distributions. Huang et al. \cite{huang2023most} prove the equivalence between the Onsager-Machlup action functional and the corresponding Markovian bridge process. They  refer 
\begin{align}
    \mbox{d}Y_t=[f(Y_t)+\sigma_t^2 \nabla \log p(x_T,T|Y_t,t)]\mbox{d}t+\sigma_t \mbox{d}W_t,\,\, t\in [0,T] \label{bridge sde}
\end{align}
as the bridge SDE associated with the system (\ref{prior sde}). This equation is originally derived by Doob \cite{doob1957conditional} from the probabilistic perspective via the Doob \textit{h-transform} of SDE (\ref{prior sde}). The relationship between the most probable transition pathway and Markovian bridge is shown in \cite{huang2023most} (Theorem 3.9) as follows

\newtheorem{theorem}{Theorem}[section]
\begin{theorem}
    Suppose that the assumptions $\mathcal{H}.1$ and $\mathcal{H}.2$ in \cite{huang2023most} hold.\\
    (\textit{i}) The most probable transition path of the system (\ref{prior sde}) coincides with the most probable transition path of the associated bridge SDE (\ref{bridge sde}).\\
    (\textit{ii}) A path $\psi^*\in C^2_{x_0,x_T}[0,T]$ is the most probable transition path of (\ref{prior sde}) if and only if it solves the following first-order ODE
    $$
\mbox{d} \psi^*(t)=\left[f\left(\psi^*(t)\right)+\sigma_t^2 \nabla \ln p\left(x_T, T \mid \psi^*(t), t\right)\right] \mbox{d} t, \quad t \in(0, T), \quad \psi^*(0)=x_0,
$$
where $p(\cdot, \cdot|\cdot, \cdot)$ is the transition density of the solution process of (\ref{prior sde}), and $C^2_{x_0,x_T}[0,T]$ denotes the space of all continuous $\mathbb{R}^n$-valued curves on time interval $[0, T]$ connecting $x_0$ and $x_T$.
\end{theorem}

\section{Tipping indicator in probability space}

Sudden transitions to a different state, once a threshold is passed, is known as tipping.
The tipping point refers to the critical point in a situation, process, or system beyond which a significant and often unstoppable effect or change takes place. Tipping points commonly signify instances where a system undergoes abrupt, swift, and occasionally irreversible transformations.
Mathematically, a tipping point is regarded as a critical state connecting the states before and after the bifurcation in a complex dynamical system. However, the general concept of the tipping point is not limited to transitions through bifurcation points; it can also result from abrupt transitions caused by instantaneous internal structural switches or parameter alterations. Ashwin et al. \cite{ashwin2012tipping} identify three mechanisms of tipping: bifurcation-induced, noise-induced,  and rate-induced tipping. 

Warning signals for critical transitions hold significant importance due to the considerable challenges involved in restoring a system to its previous state post such transitions \cite{folke2004regime, feng2024early, zhang2024early}. Various early warning indicators could be formulated with methods in statistics, bifurcation theory, information theory, statistical physics, topology, graph theory, and so on \cite{hirota2011global, wang2012flickering,quail2015predicting}. Besides, in climate dynamics \cite{Lenton2008TippingEI}, various regions in the world have been classified as tipping elements as they are considered to be under the threat of a critical transition. The definitions of "tipping elements" here inspire us to give the following definition of critical transition and tipping indicator in probability space. 

\begin{definition}(\textbf{Tipping in probability space}).
A tipping point in probability space exits if there exists a control parameter $\rho$ with a critical value $\rho_{c}$, at which a small parameter variation ($\delta \rho > 0$)
leads to a qualitative change $C$ in our tipping indicator, i.e.,
$$|\mbox{I}(\delta \rho + \rho_{c}) - \mbox{I}(\rho_{c})| \geq C.$$
\end{definition}

This definition of action functional in probability measure space is actually easy to interpret from a geometric point of view since it has been shown that certain types of PDEs whose solution is flow of probability measures are gradient flows with respect to the Wasserstein metric. In \cite{gentil2020dynamical}, a cost function which is a perturbed version of the Wasserstein cost $\mathcal{W}_2$ is defined for any regular function on the set of probability measures. 


It is also interesting to associate Lagrangian with this. Here we denote our indicator as the action functional defined on the optimal path measure via the Schrödinger bridge. Then we get
\begin{equation}
\mathcal{W}_2\left(\rho_0, \rho_1\right)=\inf_{(\nu,\, \rho)} \int_0^1 \int_{\mathbb{R}^n}\|v(t, x)\|^2 \rho_t(d x) \mbox{d} t, \label{indicator-0}
\end{equation}
where $\nu(t,x)$ is a smooth vector field. That is the cost in fluid dynamic formulation (\ref{fluid furmulation}). Therefore, denoting $ \rho_t^*$ as the optimal path measure, we have the time-dependent early warning indicator defined as
\begin{equation}
\mbox{I}(\rho_t^{*}) = \frac{1}{t}\inf _{(\nu, \,\rho)} \int_0^t \int_{\mathbb{R}^n}\|v(t, x)\|^2 \rho_t(d x) \mbox{d} t. \label{indicator}
\end{equation}



\graphicspath{{Figures/}{logo/}}
	\begin{figure*}[!t]    
	\begin{overpic}[width=0.45\columnwidth]{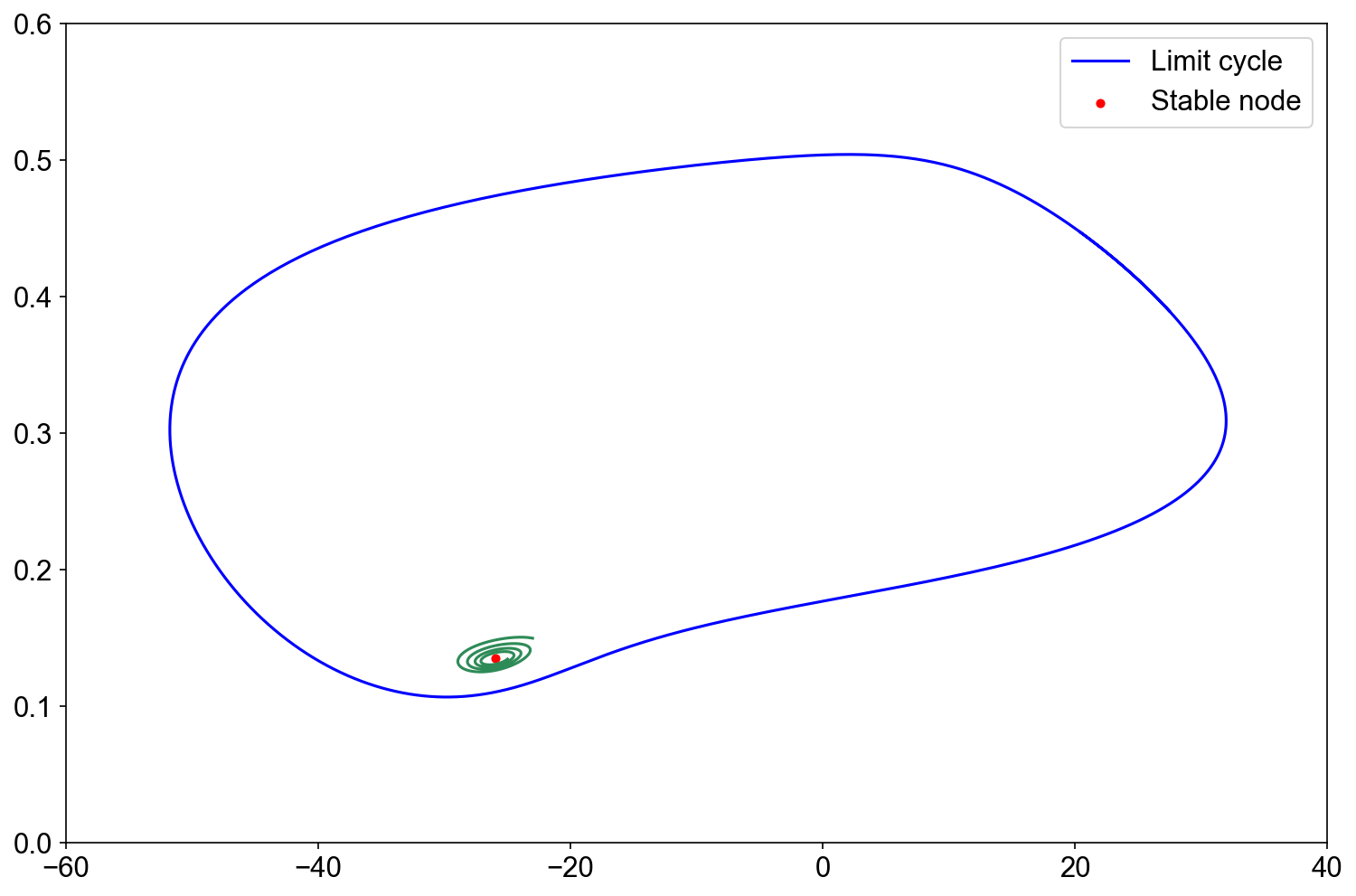}
	\put(2,70){A}
		\end{overpic}
		\hspace{2mm}
		\begin{overpic}[width=0.45\columnwidth]{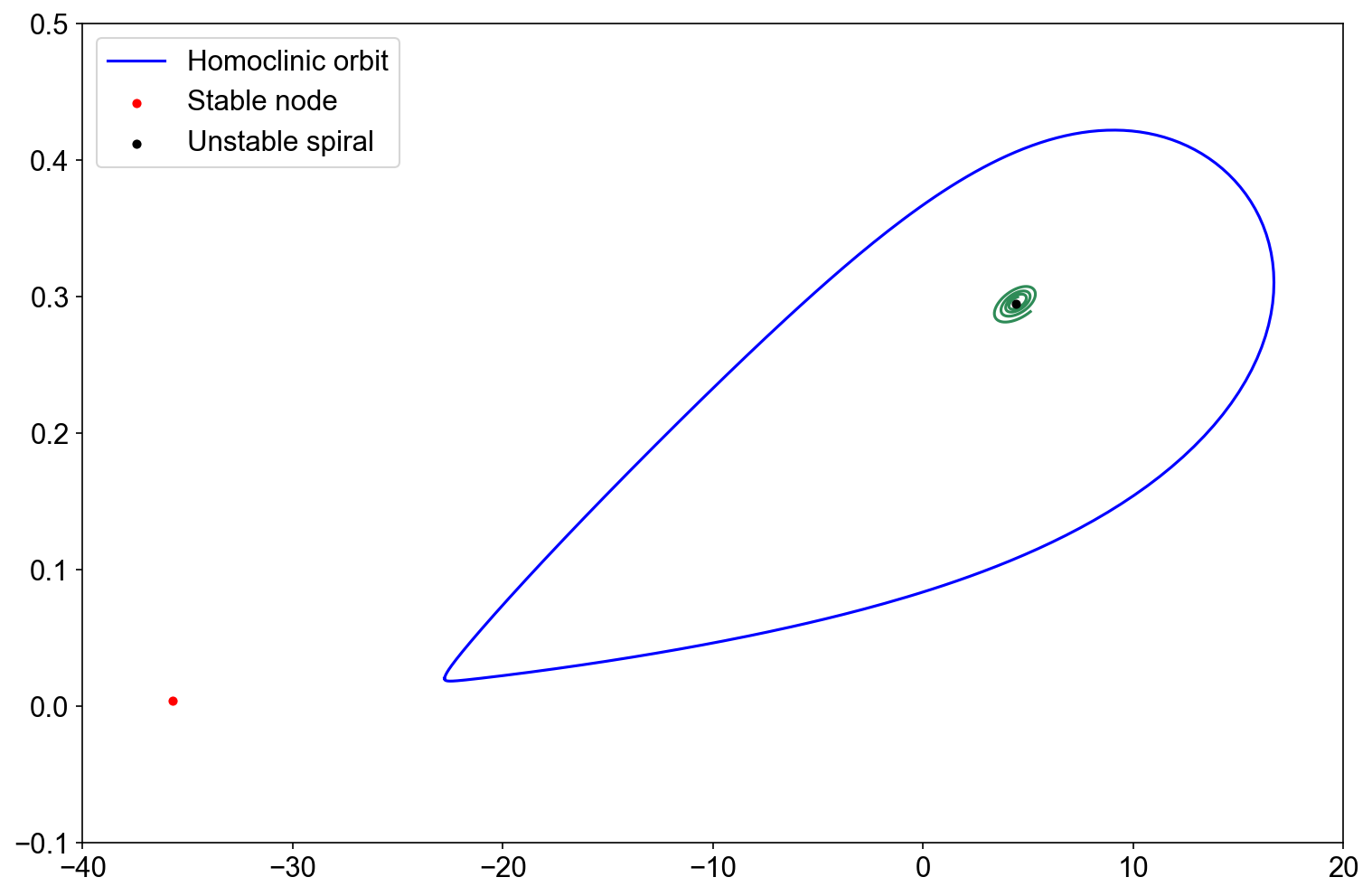}
			\put(2,70){B}
		\end{overpic}
		\captionsetup[subfigure]{labelformat=empty}
		\caption{Phase portrait of the Morris-Lecar model. (A) The stable node and stable limit cycle of the Morris-Lecar model with bifurcation parameter $I = 92 \,\,\mu A/cm^2$. (B) The stable node, unstable spiral, and stable homoclinic cycle of the Morris-Lecar model with bifurcation parameter $I = 37 \,\,\mu A/cm^2$.}
        \label{morris}
	\end{figure*}

\section{Results}

\subsection{Morris-Lecar model}\label{M-L experiement}
The Morris-Lecar model is a simplified two-dimensional biophysical model of action potential generation proposed by Catherine Morris and Harold Lecar \cite{morris1981voltage}. This model exhibits the generation of repetitive firing in certain neurons resulting from a saddle-node bifurcation on an invariant circle. The model has three channels: potassium ion channel, calcium ion channel, and leakage current channel. The model is represented as follows
\begin{numcases}{}
    & $\mbox{d} v_t =  \frac{1}{C}\left[-g_{Ca} m_{\infty}\left(v_t\right)\left(v_t-V_{Ca}\right)-g_K w_t\left(v_t-V_K\right)-g_L\left(v_t-V_L\right)+I\right] \mbox{d}t+\sigma_t \mbox{d} W_t$, \\
    & $\mbox{d} w_t =  \varphi \frac{w_{\infty}\left(v_t\right)-w_t}{\tau_w\left(v_t\right)} \mbox{d}t$,
\end{numcases}
where
\begin{align*}
    & m_{\infty}(v)=0.5\left[1+\tanh \left(\frac{v-V_1} {V_2}\right)\right], \\ 
    & w_{\infty}(v)=0.5\left[1+\tanh \left(\frac{v-V_3}{V_4}\right)\right], \\
    & \tau_w(v)=\left[\cosh \left(\frac{v-V_3}{2 V_4}\right)\right]^{-1}.
\end{align*}
Here, $W_t$ is a Brownian motion, and $\sigma_t$ is the nonnegative noise intensity. The variable $v_t$ represents the membrane potential, and $w_t \in [0,1]$ represents the activation variable of potassium ion, reflecting the evolution of the opening probability of potassium ion channels. The parameters $g_{Ca}$, $g_K$, and $g_L$ are the maximum conductances of the calcium channel, potassium channel and leakage current channel, respectively. The $V_{Ca}$, $V_K$ and $V_L$ are the reversal potentials of the above channels, respectively. The parameter $I$ represents the impressed input current from the environment, and $m_\infty$ and $w_\infty$ are the steady-state values of the opening probability ${Ca}^{2+}$ channel and $K^+$ channel, respectively. Parameter $C$ is the membrane capacitance, and $\phi$ represents the change between the fast and slow scales of the neuron. $V_1$, $V_2$, $V_3$, and $V_4$ are the parameters selected for the clamp data.



\begin{figure*}[htpb]
\begin{overpic}[width=0.47\columnwidth]{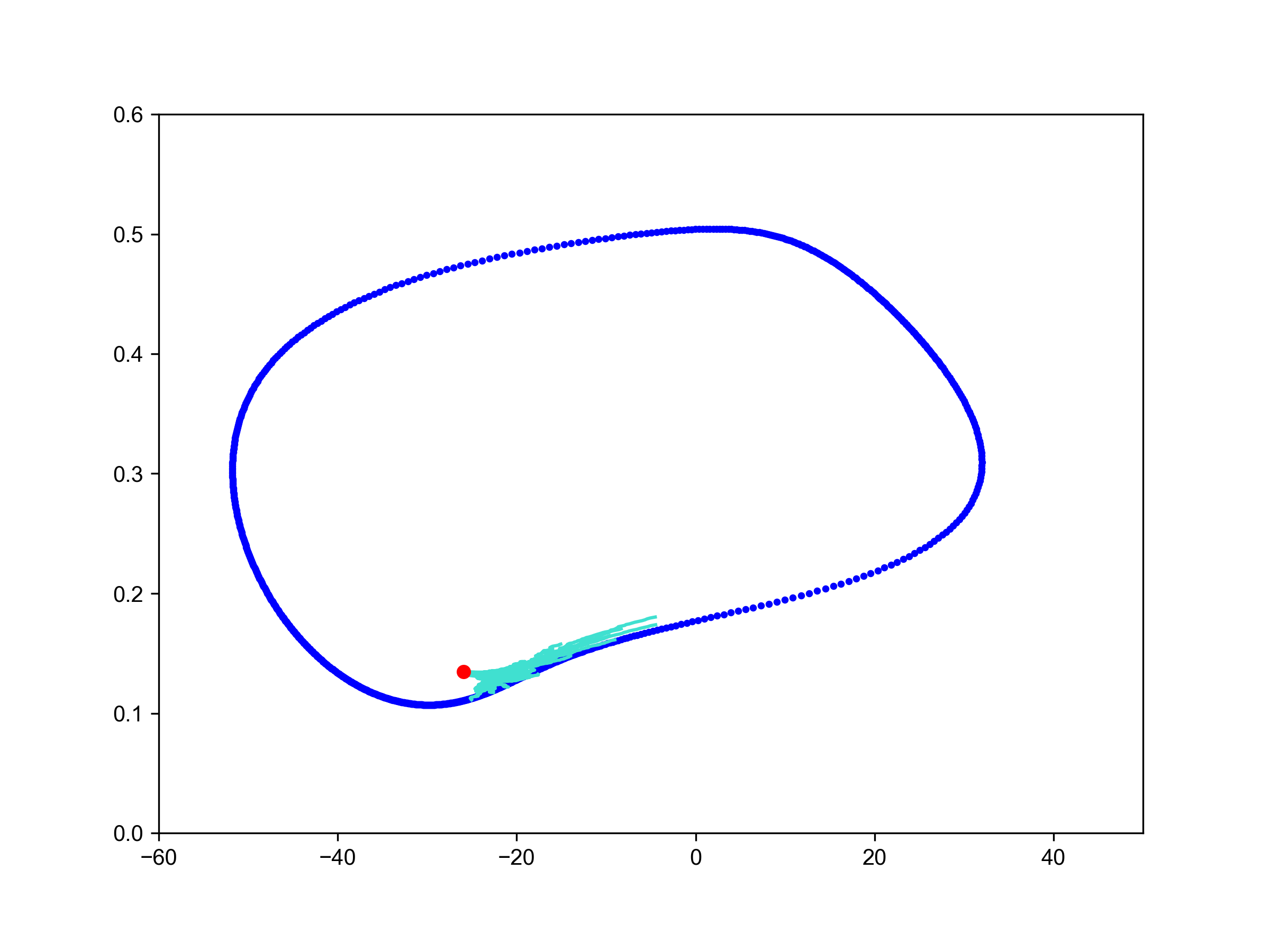}
	\put(10,73){A}
		\end{overpic}
		\hspace{2mm}
		\begin{overpic}[width=0.47\columnwidth]{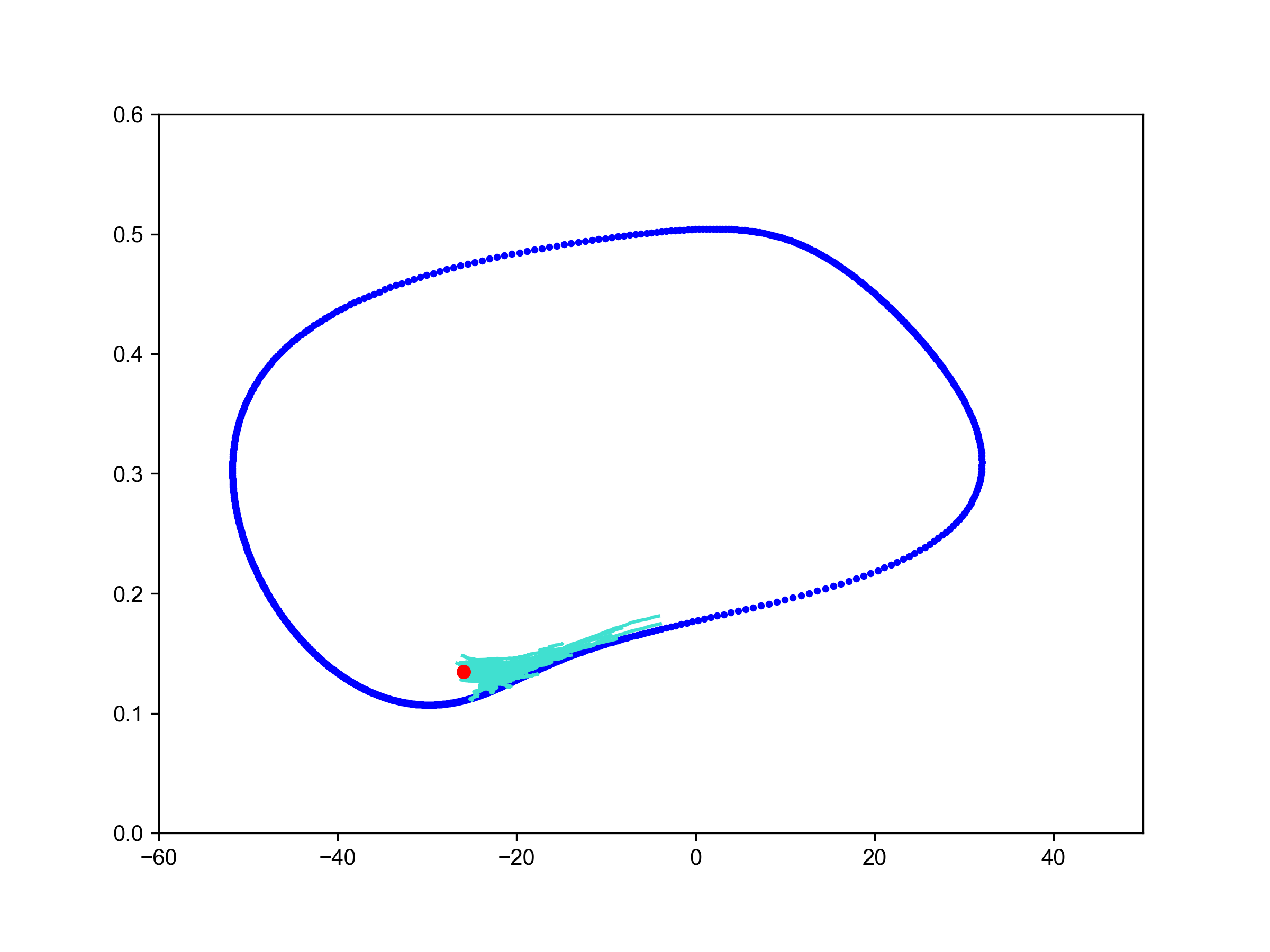}
			\put(10,73){B}
		\end{overpic}

  \begin{overpic}[width=0.47\columnwidth]{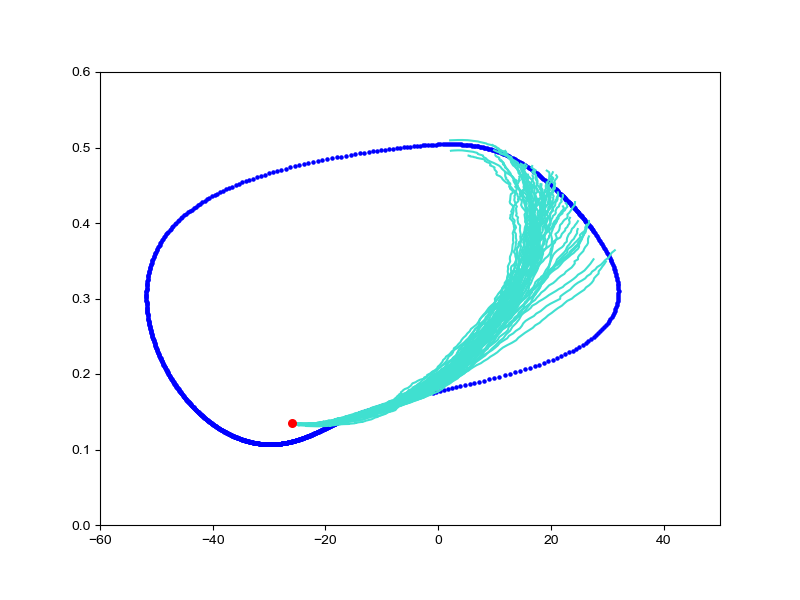}
			\put(10,73){C}
		\end{overpic}
		\hspace{2mm}
		\begin{overpic}[width=0.47\columnwidth]{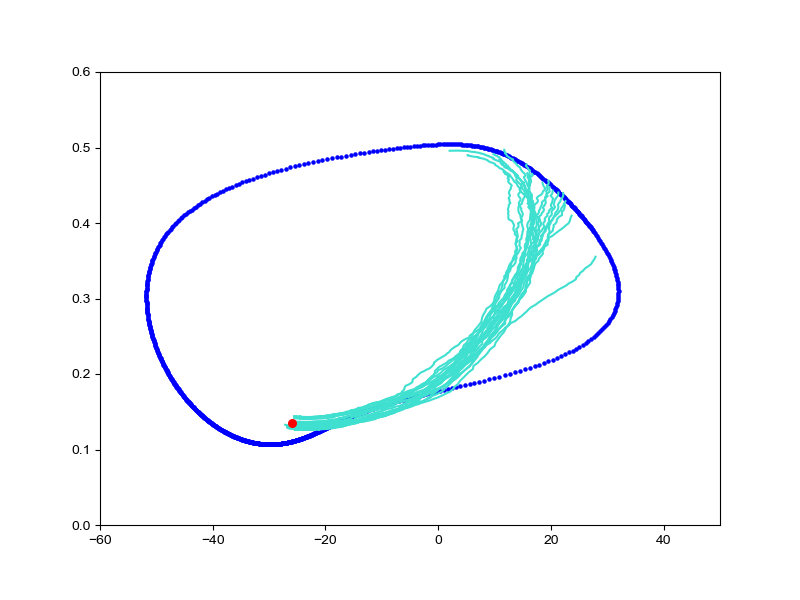}
			\put(10,73){D}
		\end{overpic}
		\captionsetup[subfigure]{labelformat=empty}
\caption{Transition path dynamics between a stable state and invariant manifold in Morris-Lecar system. Upper row: evolutionary pathways when terminal time T = 20 and noise strength g = 0.3, from forward SDE (A) and backward SDE (B).
Bottom row: evolutionary pathways when time T = 20 and noise strength g = 0.5, from forward SDE (C) and evolutionary pathways from backward SDE (D).}
\label{paths in Morris-Lecar}
\end{figure*}



Figure~\ref{morris}\textcolor{blue}{A} shows two invariant sets of the Morris-Lecar model under Class $\mbox{I}$ parameters, which is a bistable system for some parameters, with a unique stable equilibrium and a stable limit cycle. In addition, it is worth noting that there are also unstable periodic solutions (the small limiting period in Figure~\ref{morris}\textcolor{blue}{A}). The trajectory separates these initial conditions approaching the stable fixed point from those approaching the large stable limit cycle. The stable limit cycle corresponds to the continuous oscillation state of the neuronal system described by the Morris-Lecar model, while the unique fixed point corresponds to the resting state of the neuronal system. Figure~\ref{morris}\textcolor{blue}{B} shows three invariant sets of the Morris-Lecar model, with a stable node, an unstable spiral, and a stable limit cycle under Class $\mbox{II}$ parameters. The parameter details are outlined in Table~\ref{parameter} in Appendix~\ref{class parameters}.

\begin{figure}
    \centering
    \includegraphics[width=0.95\textwidth]{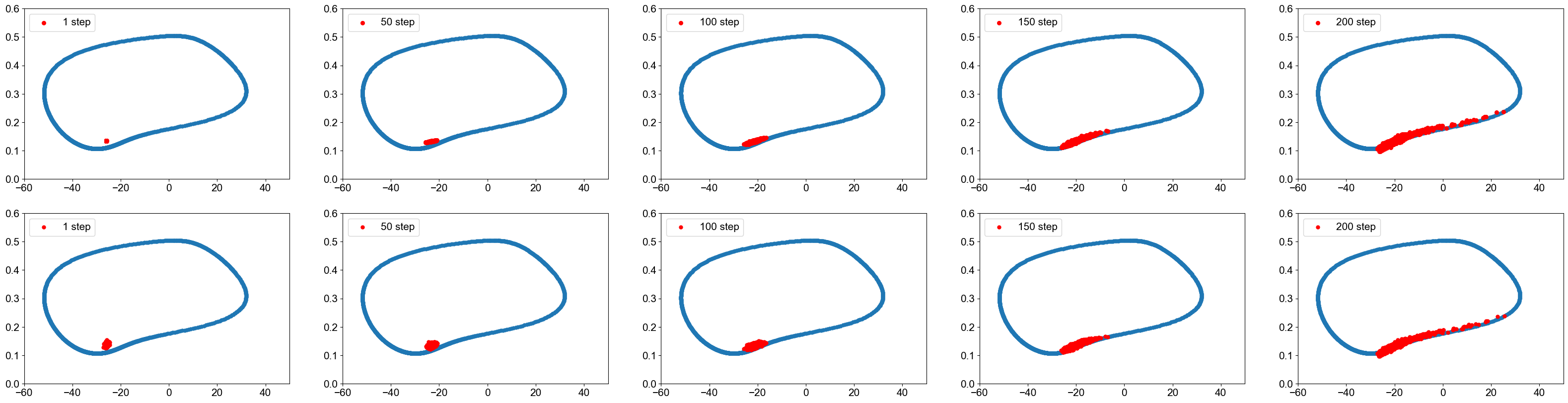}
    \caption{The evolutionary density of transition paths from the stable node to the stable limit cycle with T=20,N=200, g=0.3.}
    \label{density g=0.3}
\end{figure}

\begin{figure}
    \centering
    \includegraphics[width=0.95\textwidth]{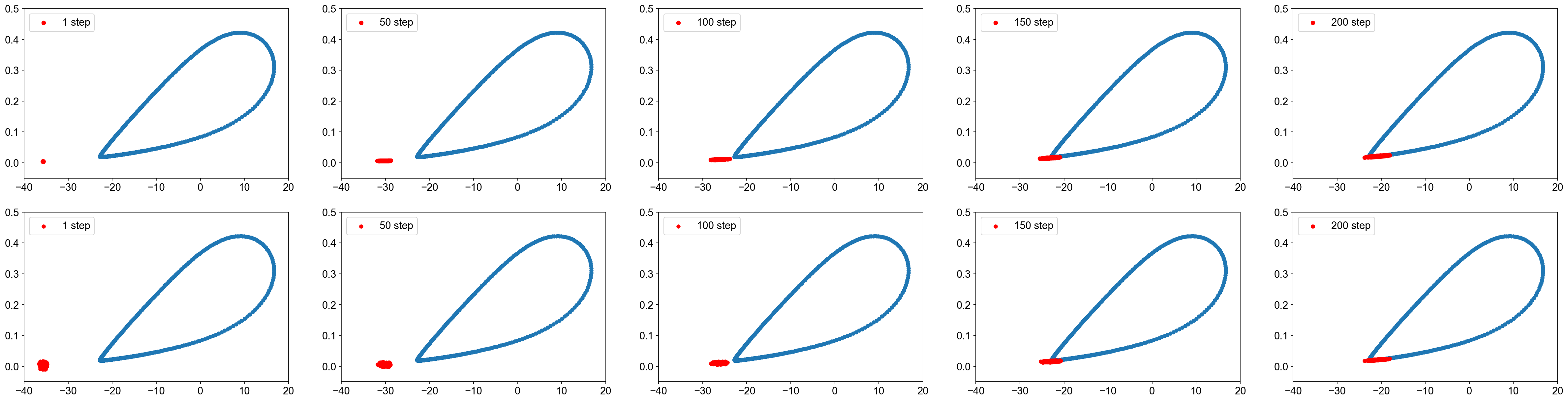}
    \caption{The evolutionary density of transition paths from the stable node to the stable homoclinic cycle with T=20, N=200,g=0.3.}
    \label{evolution-homoclinic}
\end{figure}


\begin{figure*}[!t]    
	\begin{overpic}[width=0.47\columnwidth]{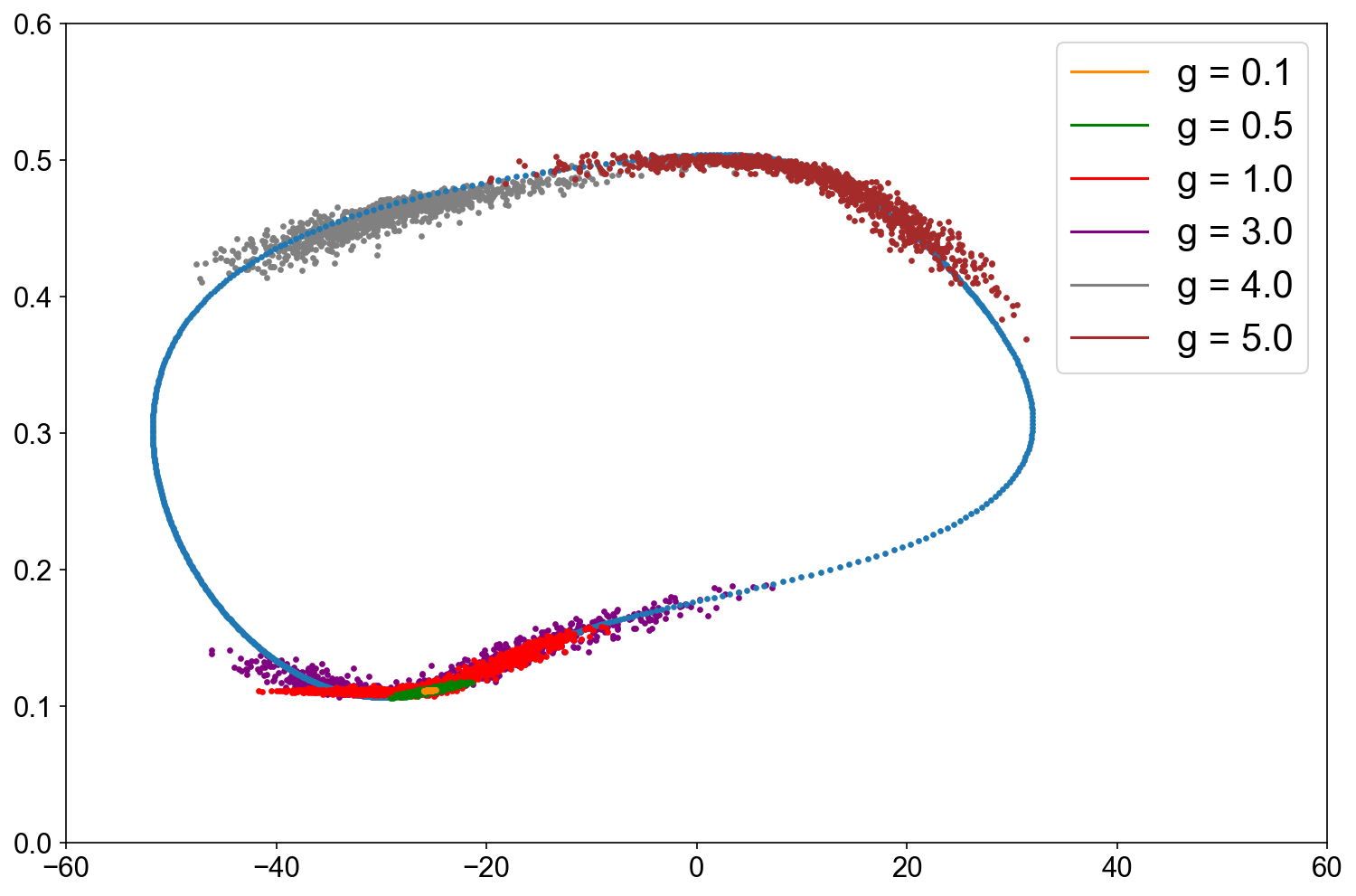}
	\put(2,70){A}
		\end{overpic}
		\hspace{2mm}
		\begin{overpic}[width=0.47\columnwidth]{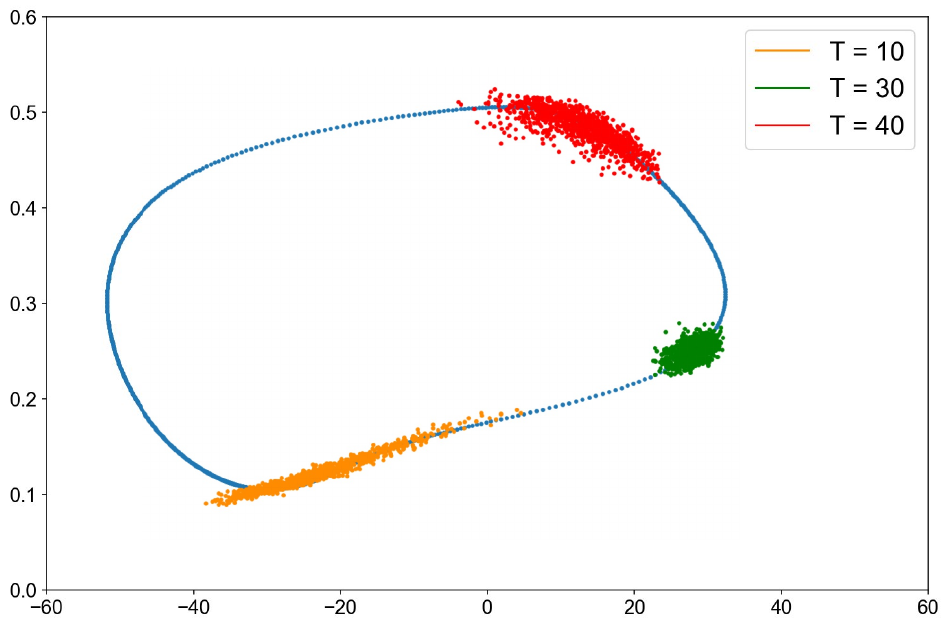}
			\put(2,70){B}
		\end{overpic}
		\captionsetup[subfigure]{labelformat=empty}
		\caption{Terminal density concentrations of the transition from the stable node to the stable limit cycle, with (A) $\mathrm{T} = 3$ under various noise density g. (B) $\mathrm{g} = 1$ under various terminal time T.}
        \label{change}
	\end{figure*}


\begin{figure}
    \centering
    \includegraphics[width=0.95\textwidth]{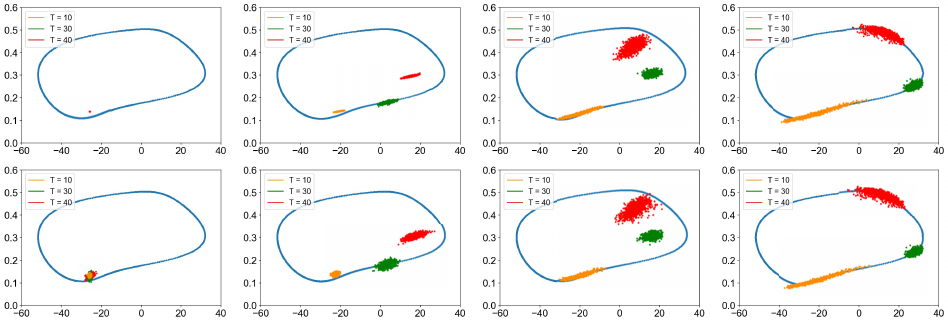}
    \caption{The evolutionary densities of transition paths dynamics over time for various terminal time $\mathrm{T} = 10, 30$ and $40$. }
    \label{evolution-vary-time}
\end{figure}


In a theoretical framework, we could endeavor to simulate the transition dynamics by aiming for an invariant measure on the stable limit cycle as the terminal time approaches infinity. Achieving precise simulation outcomes necessitates sufficiently small time steps and a large terminal time, ensuring that accumulated errors remain within acceptable bounds.    Failure to adhere to these requirements may result in the approximation of the target set being limited to a subset of the entire limit cycle. 
However, the exigency of retaining all intermediate results incurs a considerable computational cost, leading to memory constraints. Consequently, given the computational resources, we confine the target set to a subset of the limit cycle.

Figure~\ref{paths in Morris-Lecar} shows the transition path dynamics between the stable equilibrium and the invariant manifold, i.e. the limit cycle. We can see that as the noise increases, the transition paths will reach the target area at further distances. We can also see the terminal time effect of different T values. Figure~\ref{density g=0.3} and Figure~\ref{evolution-homoclinic} show the time evolutionary density of transition paths from the stable node to the stable limit cycle or the stable homoclinic cycle. The consistency of evolutionary density patterns between the forward and backward pathways means the convergence of our numerical solution to the Schrödinger bridge problem.


 


\begin{figure*}[htbp]    
	\begin{overpic}[width=0.47\columnwidth]{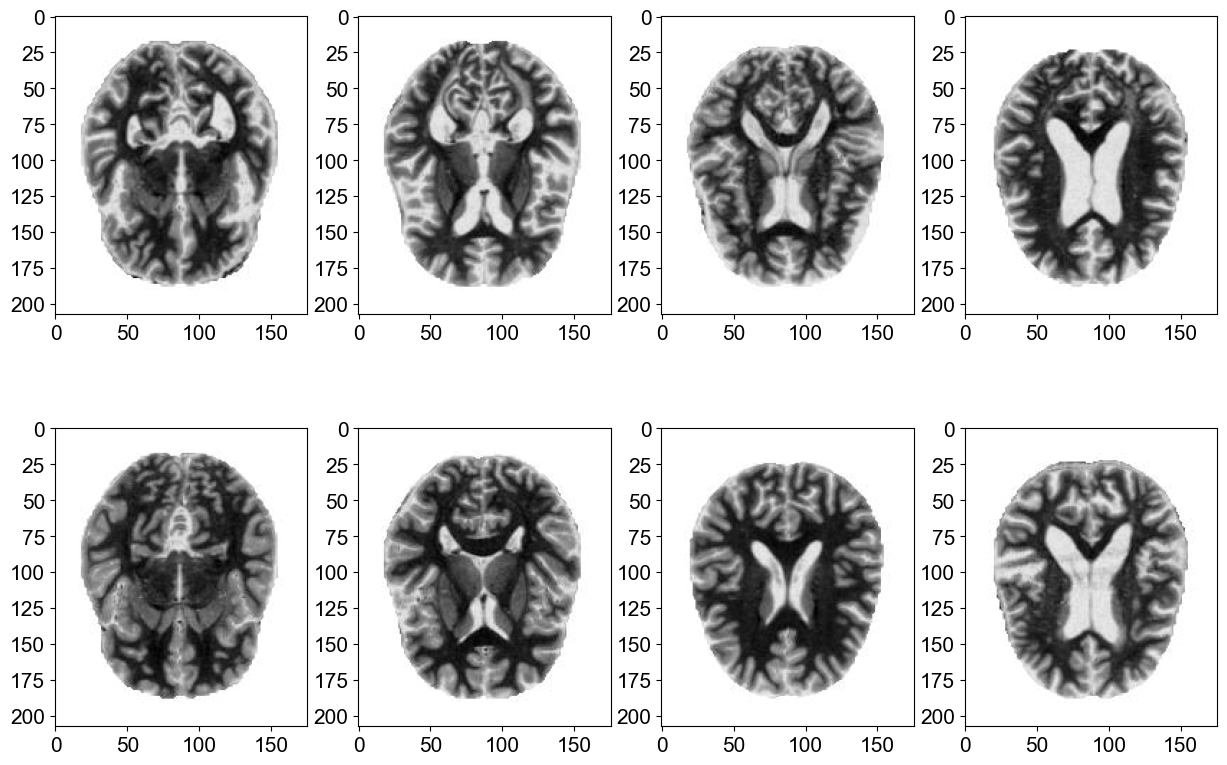}
	\put(2,70){A}
		\end{overpic}
		\hspace{2mm}
		\begin{overpic}[width=0.45\columnwidth]{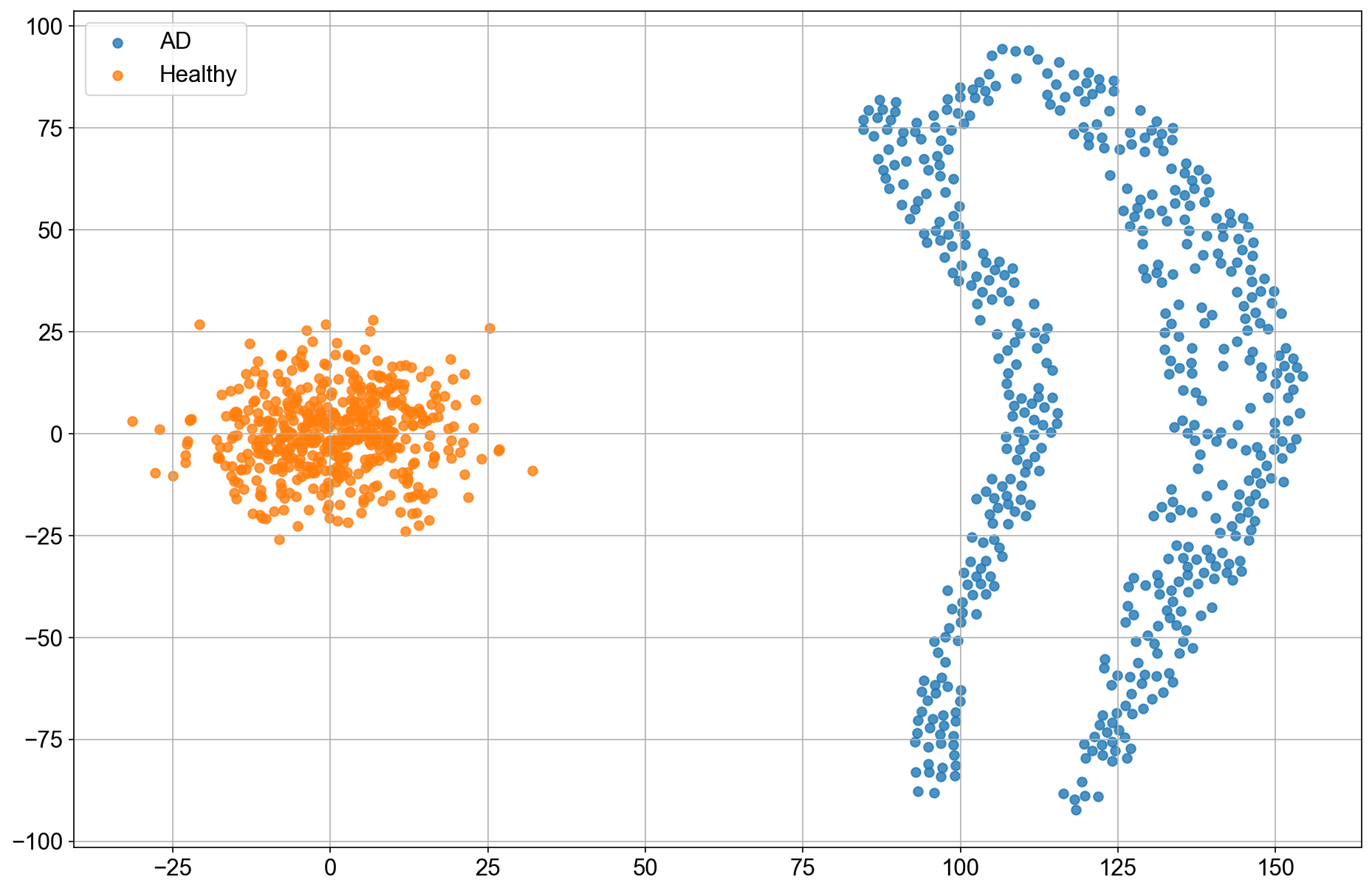}
			\put(2,73){B}
		\end{overpic}
		\captionsetup[subfigure]{labelformat=empty}
		\caption{The figure contrasts ADNI data: (A) mild disease and healthy brain samples, and (B) T-SNE visualization of healthy (orange) vs. mild disease (blue) samples.}
        \label{ADNI}
	\end{figure*}


Figure~\ref{change}\textcolor{blue}{A} and Figure~\ref{change}\textcolor{blue}{B} show the terminal density concentrations of the transition from the stable node to the stable limit cycle with various noise and terminal time. Figure~\ref{evolution-vary-time} presents the evolutionary densities of transition paths dynamics over time for various terminal time $\mathrm{T} = 10, 30$ and $40$. 
For more detail effects of the noise strength, please see Figure~\ref{different g} in Appendix~\ref{change g transition}.






\begin{algorithm}[t]
\caption{Training process for ADNI data}\label{semi ot}
\KwIn{boundary distributions $p_{\mbox{data}}$ and $p_{\mbox{prior}}=\sum_{i=1}^n \nu_i \delta(y-y_i)$, parameterized policies $Z(\cdot,\cdot;\theta)$ and $\hat{Z}(\cdot,\cdot;\phi)$,  number of Monte Carlo samples $N$, positive integer $s$.}
Initialize $h = (h_1,h_2,\dotsc,h_n) \gets (0,0,\dotsc,0).$

    \While{\textnormal{$E(h)$ has not decreased for $s$ steps}}
    {   
     Generate $N$ uniformly distributed samples from$\{x_j\}_{j=1}^N$.\\
     Calculate $\nabla h = (\hat{w}_i(h) - \nu_i)^T.$\\
     $\nabla h = \nabla h - mean(\nabla h)$.\\
     Update $h$.
    }

\For{$k = 1$ \KwTo{$K$}}
    {   
     Sample $X_{t\in [0,T]}$ from (\ref{forw}), where $x_0\sim p_{\mbox{data}}$ (computational graph retained)\;
     compute $\mathcal{L}$ with \eqref{new loss}\;
     Update($\theta, \phi$) with $\nabla_{\theta,\phi}\mathcal{L}_{SB}(x_0;\theta,\phi)$
    }
\end{algorithm}

\begin{figure*}[htb]
	\begin{overpic}[width=0.45\columnwidth]{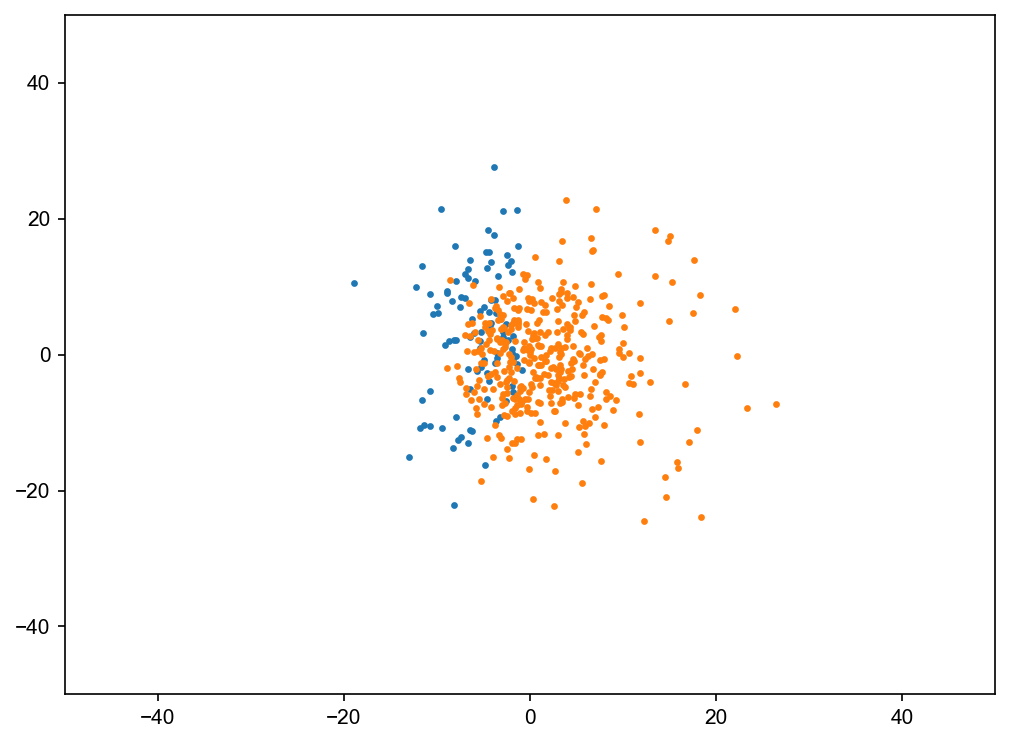}
	\put(5,80){A}
		\end{overpic}
		\hspace{2mm}
		\begin{overpic}[width=0.47\columnwidth]{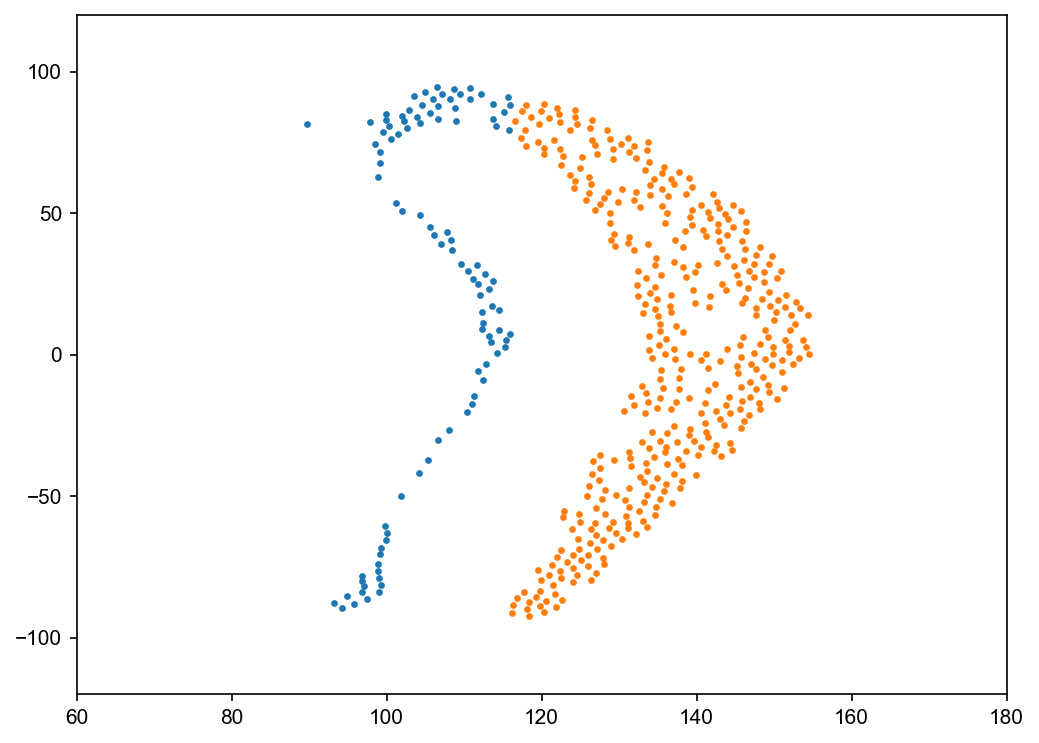}
		\put(6,77){B}
		\end{overpic}
		\captionsetup[subfigure]{labelformat=empty}
		\caption{The corresponding data points (with the same color) between initial sets (A) and target sets (B) via cell decomposition.}
        \label{part}
	\end{figure*}

\begin{figure}[htb]
    \centering
    \includegraphics[width=1.0\textwidth]{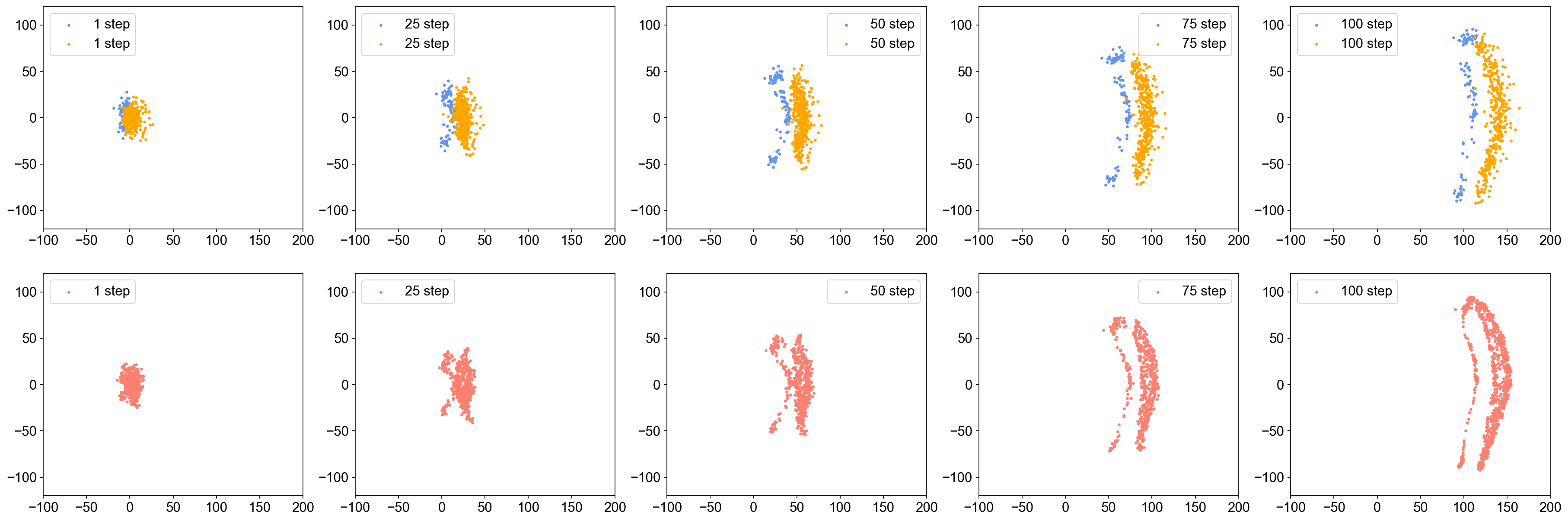}
    \caption{The evolutionary density from the healthy state to early Alzheimer's state.}
    \label{transition from normal to mild}
\end{figure}

\subsection{Alzheimer's disease neuroimage}


AD is distinguished by a progressive and seemingly continuous pathological progression, marked by the development of neurodegeneration exhibiting distinct spatial and temporal dynamics. Specific characteristics can strongly indicate the diagnosis of AD.  In AD, the brain typically exhibits at least moderate cortical atrophy, particularly pronounced in multimodal association cortices and limbic lobe structures. Structural neuroimaging with MRI can be employed to evaluate atrophy as a metric of neurodegeneration in AD. We conduct our experiments on dataset originating from \href{https://www.kaggle.com/datasets/kaushalsethia/alzheimers-adni}{ADNI} (Samples in Figure \ref{ADNI}\textcolor{blue}{A}).

Researchers \cite{masquelier2023new} have determined that the presence of abnormally tangled tau proteins serves as a sign of numerous instances of Alzheimer's disease. They noted that the tau protein undergoes a crucial \textit{tipping point} transition, shifting from a healthy to a pathological state. Upon crossing this tipping point, tangles are rapidly formed. Once this threshold is surpassed, the brain transitions into the irreversible clinical phase, losing opportunities for intervention \cite{simons2023tipping}. Therefore, the detection of the tipping point in Alzheimer's disease holds significant importance, enabling medical professionals to intervene effectively before the onset of the clinical phase.

Due to the difficulty of detecting the tipping point in high-dimensional time series, we need to reduce the dimension of brain MRI data to the plane. Firstly, the Variational AutoEncoder (VAE) \cite{kingma2013auto} model is employed to extract features from the representation space, reducing the image data to 20-dimensional vectors. Subsequently, these 20-dimensional vectors are projected into a 2-dimensional space utilizing the t-distributed Stochastic Neighbor Embedding (t-SNE) \cite{van2008visualizing} method. The following is the visualization of healthy and AD patients after VAE and t-SNE reduction (Figure \ref{ADNI}\textcolor{blue}{B}).

\begin{figure*}[htb]
	\begin{overpic}[width=0.47\columnwidth]{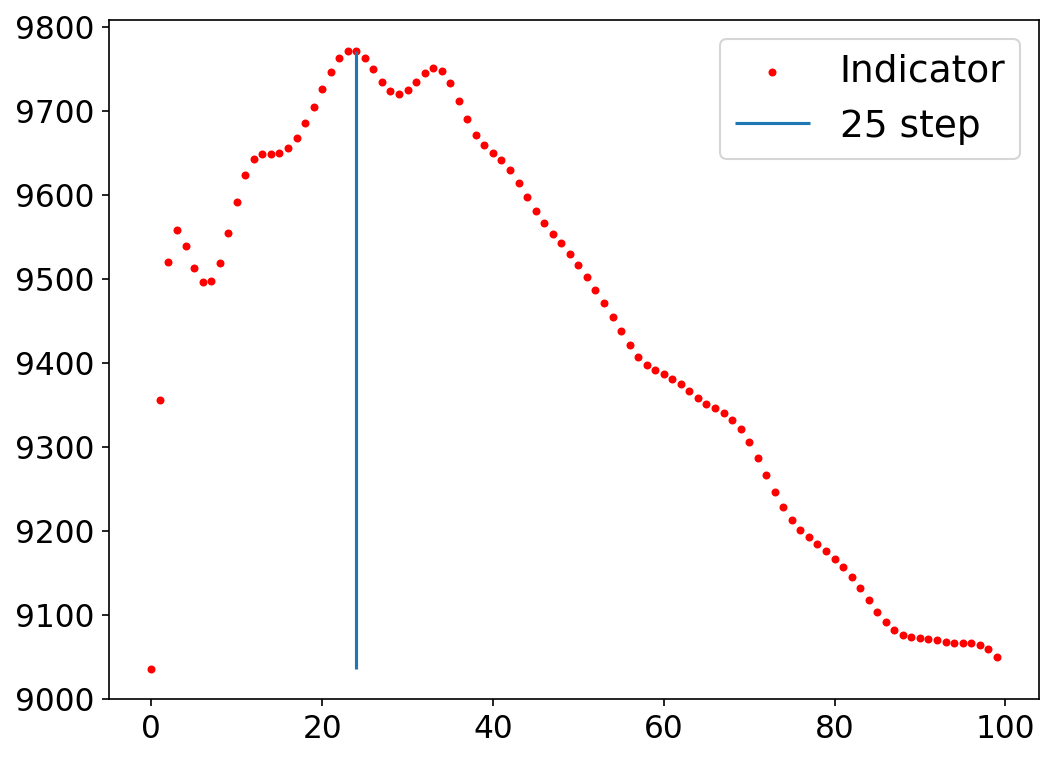}
	\put(3,78){A}
		\end{overpic}
		\hspace{2mm}
		\begin{overpic}[width=0.47\columnwidth]{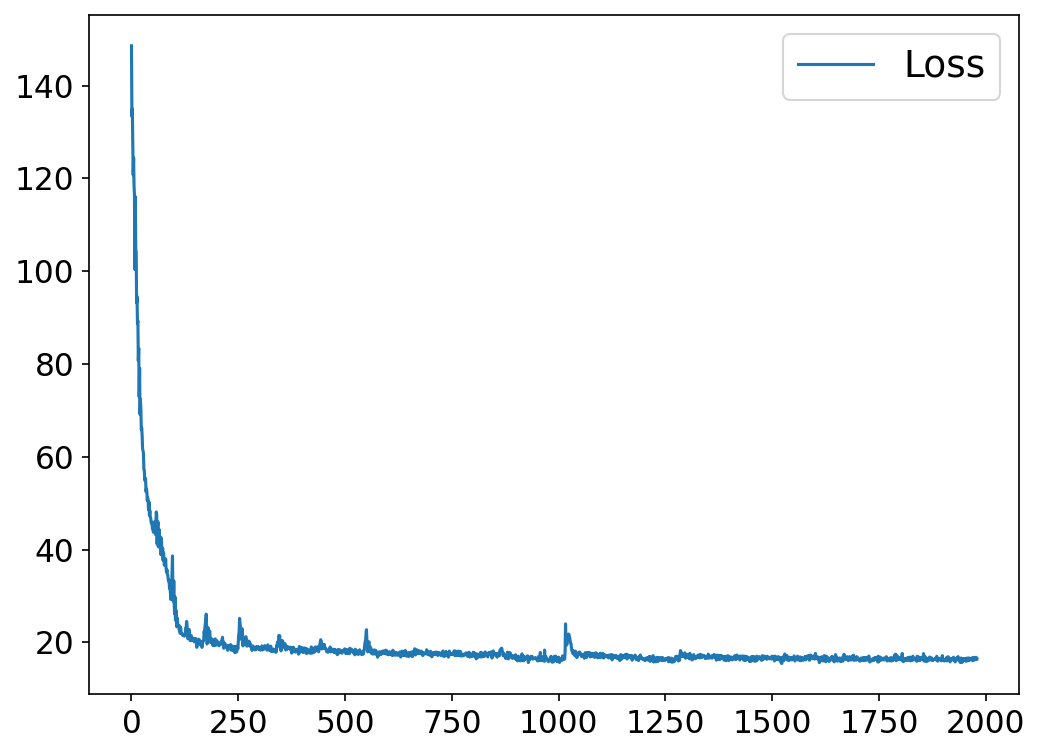}
		\put(3,78){B}
		\end{overpic}
		\captionsetup[subfigure]{labelformat=empty}
		\caption{The action functional indicator defined in Eq.(\ref{indicator}) with respect to time (A) and the loss function over iterations (B).}
        \label{cost and loss}
	\end{figure*}



By observing Figure \ref{ADNI}\textcolor{blue}{B}, it becomes evident that the data pertaining to early-stage mild cognitive impairment exhibits a non-convex configuration post dimensionality reduction. 
According to Brenier’s Theorem (\cite{brenier1987polar,brenier1991polar}), any transport map can be decomposed into a measure preserving map and a solution to the Monge-Amp$\grave{\mathrm{e}}$re equation. Therefore,  under the $L^2$ cost function, the regularity of the solution to the Monge-Amp$\grave{\mathrm{e}}$re equation determines the continuity of the transport map.
If the support set characterizing the target distribution lacks connectivity or convexity, the optimal transport map will manifest discontinuous property. Nonetheless, conventional neural networks are confined to approximate continuous mappings, thus engendering inherent conflicts that precipitate convergence challenges and the phenomenon of pattern collapse. To redress this quandary, we embrace the methodology delineated in \cite{an2019ae,lei2019mode}, predicated upon the discrete Brenier theory to approximate the continuous Brenier potential.

Specifically, denote the initial measure as $\mu$ over the convex compact domain $\Omega$, and the objective measure $\nu$ as an empirical measure, i.e. $\nu = \sum_{i=1}^n\nu_i\delta(y-y_i)$, where $Y = \{y_1, y_2, \cdots, y_n \} $, and $\nu_i $ meet $\sum_{i=1}^n\nu_i = \mu(\Omega)$. Each training sample $y_i$ corresponds to the support plane of the Brenier potential energy function, denoted by
$$
\pi_{h,i}(x) := \langle x,y_i \rangle + h_i,
$$
and the height vector is defined as $h = (h_1,h_2,\cdots, h_n)$. Denote $\omega_i(h)= \mu(W_i(h))=\int_{W_i(h)\cap \Omega} \mbox{d}\mu$, by semi-discrete Optimal Transport Map \cite{lei2019mode}, the vector $h$ is the unique minimizer of the energy $E(h)=\int_0^h \sum_{i=1}^n \omega_i(\eta)\mbox{d}\eta_i-\sum_{i=1}^n h_i\nu_i$.

\begin{figure*}[htb]
\centering
	\begin{overpic}[width=0.9\columnwidth]{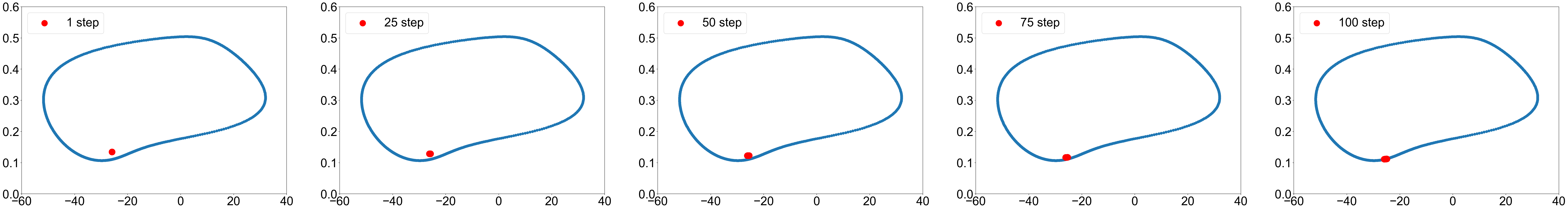}
	\put(0.5,14){A}
		\end{overpic}
  \\
  \hspace{4cm}
		\begin{overpic}[width=0.9\columnwidth]{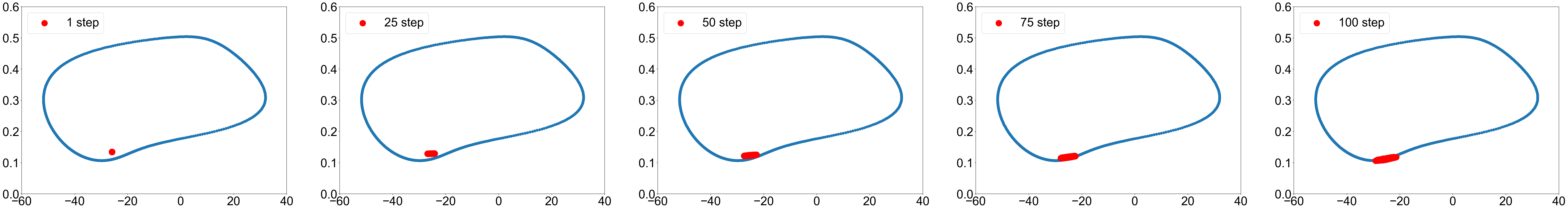}
		\put(0.5,14){B}
		\end{overpic}
\\
  \hspace{4cm}
        \begin{overpic}[width=0.9\columnwidth]{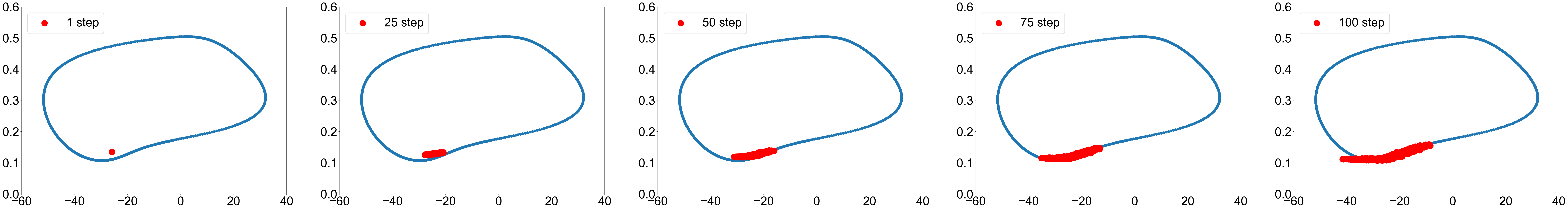}
		\put(0.5,14){C}
		\end{overpic}
\\
  \hspace{4cm}
   \begin{overpic}[width=0.9\columnwidth]{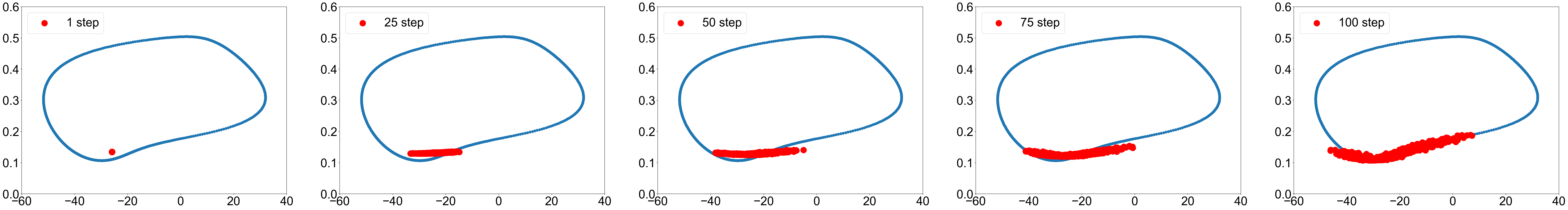}
		\put(0.5,14){D}
		\end{overpic}
\\
  \hspace{4cm}
\begin{overpic}[width=0.9\columnwidth]{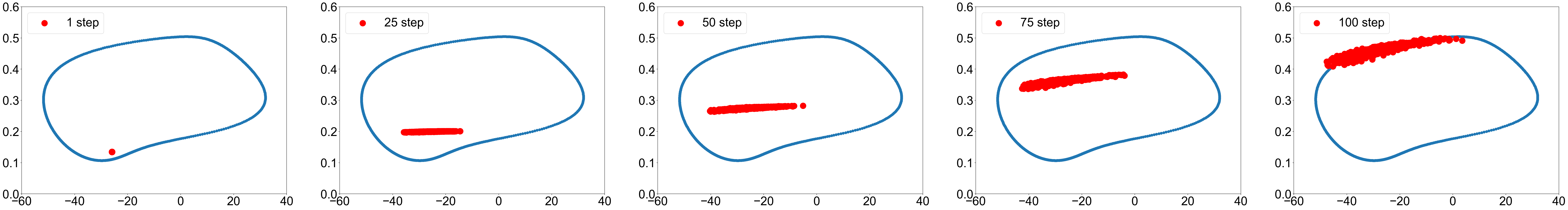}
		\put(0.5,14){E}
		\end{overpic}
\\
  \hspace{4cm}
  \begin{overpic}[width=0.9\columnwidth]{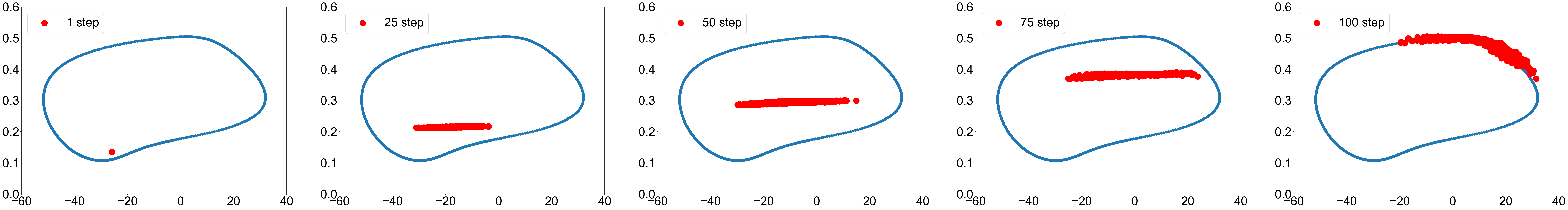}
		\put(0.5,14){F}
		\end{overpic}
		\captionsetup[subfigure]{labelformat=empty}
		\caption{The evolutionary densities of transition paths dynamics over time for various noise density g. (A) $g=0.1$. (B)$g=0.5$. (C)$g=1$. (D)$g=3$. (E)$g=4$. (F)$g=5$.}
\label{different g}
	\end{figure*}







The figure of the Brenier potential function is a convex polyhedron with each support plane $\pi_{h,i}(x)$ corresponding to one of the faces of the polyhedron. Under a semi-discrete transport map $T:\Omega\rightarrow Y$, a cell decomposition is induced $\Omega=\cup_{i=1}^n W_i$, such that
every $x$ in each cell $W_i$ is mapped to the target $y_i$, $T:x\in W_i \rightarrow y_i$. The map $T$ is measure
preserving, denoted as $T_\# \mu = \nu$, if the transport map satisfies $T(W_i)=y_i$ and $\mu(W_i)=\nu_i$. The specific training process is shown in algorithm \ref{semi ot}.

By McCann interpolation, the target distribution satisfies $\rho_1(x_1)=(\nabla \varphi)_\# \rho_0(x_0)$, where $\varphi$ is the Kantorovich potential. Under $L^2$ cost function, the Brenier theory indicates that the optimal transport $T$ and Kantorovich potential $\varphi$ meets $T(x)=Id(x)-\nabla \varphi(x)$. Furthermore, we obtain $T(x)=\nabla (\pi_{h,i}(x))=y_i$. Therefore, the $\rho_1(x_1)=(\nabla \varphi)_\# \rho_0(x_0)=(Id - T)_\# \rho_0(x_0)$. Then we modify the loss function (\ref{loss function}) into
\begin{equation}
\mathcal{L} = - \ln{p_0(x_0)}= \int_0^1 \mathbb{E}\big[\frac{1}{2}(\|\sigma_t \nabla\varphi_t\| +\|\sigma_t \nabla\widehat{\varphi}_t\|)^2 + \nabla_x\cdot (\sigma_t^2 \nabla\widehat{\varphi}_t-f)\big]\mbox{d}t - \mathbb{E}[(X_1-y_i)^2], \label{new loss}
\end{equation}
where $x_0\in W_i$, and $X_1$ is the propagated terminal state which starts from $x_0$ by forward SDE (\ref{forw}).

Due to the non-convex geometric shape of target data, we divide it into two connected regions, as the blue part and the orange part in Figure \ref{part}\textcolor{blue}{B}. Then by Semi-discrete Optimal Transport Map, we obtain the corresponding two parts in initial data, which is shown in Figure \ref{part}\textcolor{blue}{A}. Next, we propagate the FBSDEs for the two parts respectively. The algorithm is shown in Algorithm \ref{semi ot}.

We employed reduced MRI images from 1100 healthy individuals as the initial distribution and from 750 AD patients as the target distribution. The pathway dynamics from the initial to the target distribution are established via the Schrödinger bridge theory, which could be viewed as an entropy-regularized optimal transport problem. To identify the tipping point from the healthy brain state to the AD brain state, we compute the estimates of $\nabla \ln{\phi_t}$ and density $\rho_t$ via the loss function described in Eq.(\ref{loss function}).

The embedding densities of the two datasets can now be employed to compare tipping phenomena. We choose terminal time T = 1, N = 100 and the noise intensity g is adaptive to time change\cite{nichol2021improved}. From Figure~\ref{transition from normal to mild}, we can see the time evolutionary density from the healthy state to the Alzheimer's state. The consistent evolution of the first and second rows across steps in Figure~\ref{transition from normal to mild} represents that the determined density involution converges, a theoretical guarantee provided by \cite{leonard2013survey,albergo2023stochastic}. It's obvious that the geometric of samples undergo a sudden change around the 25-th step, aligning with the abrupt tipping indicated by our defined indicator in (\ref{indicator}). That implies that a critical transition in the brain state occurs around 25-th step, emphasizing the need for intervention prior to the tipping point.

Figure~\ref{cost and loss}\textcolor{blue}{A} illustrates the temporal consistency of the indicator's maximal value.with the most probable splitting at 25-th step. This indicator, defined in Eg. (29), alignswell with the entropy production rate (EPR), indicating that the system experiences themost significant disorder around 25-th step. The optimal transition path between thetwo sets via the Schrodinger bridge, is guaranteed to be convergent based on the forward-backward SDE with the loss function in Figure~\ref{cost and loss}\textcolor{blue}{B}, as well as the forward and backward consistency in Figure~\ref{transition from normal to mild}.




\begin{table*}[htbp]
\centering
\caption{The parameters for Morris-Lecar model.}
\label{parameter}
\begin{tabular}{ccc}
   \toprule
   Parameter & Class $\mbox{I}$ (Hopf) & Class $\mbox{II}$ (Homoclinic) \\
   \midrule   
   $g_{Ca}$ & $4.4\,\,$ & $4$\\
   $g_K$ & $8\,\,$ & $8$\\
   $g_L$ & $2\,\,$ & $2$\\
   $E_{Ca}$ & $120\,\,$ & $120$\\
   $E_K$ & $-84\,\,$ & $-84$\\
   $E_L$ & $-60\,\,$ & $-60$\\
   $I$ & $92\,\,$ & $37$\\
   $C$ & $20\,\,$ & $20$\\
   $\phi$ & $0.04\,\,$ & $0.23$\\
   $V_1$ & $-1.2\,\,$ & $-1.2$\\
   $V_2$ & $18\,\,$ & $18$\\
   $V_3$ & $2\,\,$ & $12$\\
   $V_4$ & $30\,\,$ & $17.4$\\
   
   \bottomrule
\end{tabular}
\end{table*}

\section{Conclusion}\label{CO}
\noindent

In this paper, we introduce a novel and effective method for automatically detecting abnormalities in neural images. Our approach leverages the Schrödinger bridge with early warning indicators which facilitate us in two aspects. On the one hand, we investigate the transition path dynamics in terms of path measures between two meta-stable invariant sets. On the other hand, we study the early warning signals of real ADNI data via action functional in probability space. The ability to identify transitions between meta-stable states plays a pivotal role in predicting and controlling brain behavior and functionality. We here find that the action functional indicator in probability space is a good metric to measure the sudden discrepancy as density changes over time. 

There are still challenging problems to be solved in the future. For example, as the geometric representation of images on low-dimensional manifolds can greatly affect the convergence of our algorithm, one important issue is to develop better algorithms that take into account the geometric properties of the given image distributions. 

\section{Appendix}

\subsection{Parameters for Morris-Lecar Model}\label{class parameters}
The parameter details outlined in Section~\ref{M-L experiement} are presented in Table~\ref{parameter}. In the table, the second column corresponds to a scenario featuring one stable state and a stable cycle, while the third column corresponds to a situation with a stable node, an unstable spiral, and a stable limit cycle.

\subsection{Transition Paths Dynamics over Time}\label{change g transition}
The evolutionary densities of transition paths dynamics over time for various noise density $g$ is shown in Figure~\ref{different g}.

\section*{Author contribution}
\textbf{Peng Zhang:} Conceptualization, Methodology, Software, Formal analysis, Writing – original draft. \textbf{Ting Gao:} Methodology, Software, Formal analysis, Writing - review \& editing, Resources, Funding Support. \textbf{Jin Guo:} Methodology, Software, Visualization, Formal analysis, Writing – original draft. \textbf{Jinqiao Duan:} Project administration, Writing - review \& editing, Supervision, Funding Support.

\section*{Acknowledgements}
This work was supported by the National Key Research and Development Program of China (No. 2021ZD0201300), the National Natural Science Foundation of China (No. 12141107), the Fundamental Research Funds for the Central Universities (5003011053), the Guangdong-Dongguan Joint Research Grant 2023A151514 0016, and is partly supported by the Guangdong Provincial Key Laboratory of Mathematical and Neural Dynamical Systems. 

\section*{Conflict of interest statement}
The authors declare that they have no conflict of interest.

\section*{Ethics statement}
This article does not contain any studies with human or animal materials performed by any of the authors.

\bibliographystyle{unsrt}
\bibliography{main}

\end{document}